\documentclass[a4paper]{article}
\usepackage{amssymb}
\usepackage{amsmath}
\usepackage{amsfonts}

\newtheorem{theorem}{Theorem}[section]
\newtheorem{lemma}[theorem]{Lemma}
\newtheorem{corollary}[theorem]{Corollary}
\newtheorem{proposition}[theorem]{Proposition}

\newtheorem{remark}[theorem]{Remark}

\newtheorem{definition}[theorem]{Definition}

\newtheorem{assumption}[theorem]{Assumption}

\newenvironment{Proof}{\removelastskip\par\medskip
\noindent{\em Proof.} \rm}{\penalty-20\null\hfill$\square$\par\medbreak}

\begin{document}
\title{On the modularity of supersingular elliptic curves over certain totally real number fields}
\author{Frazer Jarvis and Jayanta Manoharmayum\thanks{The authors
thank the EPSRC for its support by means of a Research Grant from 2001--2003, when this work was done.}}

\date{{\tt a.f.jarvis@shef.ac.uk}\qquad{\tt j.manoharmayum@shef.ac.uk}}
\maketitle

\begin{large}
\end{large}

\begin{abstract}
We study generalisations to totally real fields of the methods originating with Wiles and
Taylor-Wiles (\cite{Wi},~\cite{TW}). In view of the results of Skinner-Wiles~\cite{SW2} on
elliptic curves with ordinary reduction, we focus here on the case of supersingular reduction.
Combining these, we then obtain some partial results on the modularity problem for semistable
elliptic curves, and end by giving some applications of our results, for example proving the
modularity of all semistable elliptic curves over $\mathbb{Q}(\sqrt{2})$.
\end{abstract}

\section{Introduction}
Let $E$ denote an elliptic curve over a totally real number field $F$. We say that $E$ is {\sl
modular} if there is a Hilbert modular form $f$ over $F$ of parallel weight 2 (i.e., the
corresponding automorphic representation has weight 2 at every infinite place) such that the
Galois representation associated to $E$ via its $\ell$-adic Tate module is isomorphic to an
$\ell$-adic representation associated to $f$ (see~\cite{Ca} and~\cite{Ta1}).

The approach is now standard, and originated in~\cite{Wi} and~\cite{TW}; one considers the case
$\ell=3$, uses the Langlands-Tunnell Theorem to show that the {\it reduction}
$\overline{\rho}_{E,3}$ is modular, and then proves that every (suitably constrained) lift to
characteristic 0 is modular.

Historically, the easier case has been where $\overline{\rho}_{E,3}$ is irreducible. In this
case, the deformation theory is now well understood, and this was the only case needed by Wiles
and Taylor-Wiles (\cite{Wi},~\cite{TW}). Over totally real fields, Fujiwara circulated a
manuscript~\cite{Fu} some years ago, proving an important generalisation of the method of
Taylor-Wiles, and announcing a proof of the modularity of certain elliptic curves over totally
real fields. However, there are
several hypotheses appearing in his main theorem which we hope partially to eliminate in this
work. Subsequently, Skinner and Wiles~\cite{SW2} have proven the modularity in many `nearly
ordinary' cases.

In the case where $\overline{\rho}_{E,3}$ is reducible, Skinner and Wiles~\cite{SW1} have
developed new techniques to demonstrate modularity of elliptic curves (and more general Galois
representations) over totally real fields, although these results depend on certain hypotheses on
cyclotomic extensions of $F$. Since the first version of this article was written (2002-3), Kisin
has also found stronger results (see~\cite{Ki1},~\cite{Ki2}).

\subsection*{Reduction to the semistable case}
We first remark that the modularity of all elliptic curves over totally real fields may be
reduced to proving the modularity of all semistable elliptic curves over totally real fields. The
argument is simple; by an explicit version of the semistable reduction theorem (see, for example,
\cite{taylorartin}, Lemma~2.2), an elliptic curve $E$ over a totally real field $F$ attains
semistable reduction over a finite soluble totally real Galois extension $F'/F$. (Note that
$F'/F$ will be ramified at any prime of $F$ at which $E$ has additive reduction.) The modularity
of $E_{/F}$ then follows from the modularity of $E_{/F'}$ using base-change techniques. This
argument is well-known to experts, so we omit it here.

For this reason, we restrict attention to semistable curves, and try to prove modularity. In view
of some of the applications in mind, we focus in this paper on the easiest case, where the
ramification conditions on the field are as strong as possible, but the methods should apply more
generally. Because of the results already obtained in the reducible and ordinary cases, we focus
on the supersingular case in this paper.

\subsection*{Applications}
As we are able to prove the modularity of more elliptic curves than was previously known, we can
therefore improve certain results in the literature. Following Wiles's methods~(\cite{Wi}), we
try to find fields for which we can prove modularity of all semistable curves. Wiles~(\cite{Wi},
chapter~5) uses a switch between the primes 3 and 5, which depends on the finiteness of
$X_0(15)(\mathbb{Q})$; however $X_0(15)(F)$ will generally not be finite. Other restrictions on
the field also become apparent in generalising directly his methods. However, we are able to
prove modularity of all semistable elliptic curves for the quadratic fields
$\mathbb{Q}(\sqrt{2})$ and $\mathbb{Q}(\sqrt{17})$. That we can prove such results for the first
of these fields is a piece of good fortune; the first author and Paul Meekin (\cite{JM}) have
shown that a generalisation of Fermat's Last Theorem to $\mathbb{Q}(\sqrt{2})$ would follow from
such a result. They also show that $\mathbb{Q}(\sqrt{2})$ is the only real quadratic field for
which an implication of the form `modularity implies Fermat' can be derived directly.

\subsection*{Notation}
The absolute Galois group of a field $F$ is written either as $\mathrm{Gal}(\overline{F}/F)$ or
$G_F.$ The separable algebraic closure of $F$ is denoted by $\overline{F}$. Given an extension of
fields $K\supset F$ and some representation $\rho:G_F\rightarrow \mathrm{GL}_2(*),$ we denote the
restriction of $\rho$ to the absolute Galois group of $K$ by either $\rho|_{G_K}$ or, simply, by
$\rho|_K.$ If $F$ is a number field, we denote the decomposition and inertia groups at a place
$v$ by $D_v$ and $I_v$ respectively.

Throughout, $\ell$ is an odd prime. We denote the $\ell$-adic cyclotomic character by
$\epsilon_\ell$, and its reduction, the mod $\ell$ cyclotomic character, by
$\overline{\epsilon}_\ell$. We denote by $\omega_2$ the second fundamental character of
$\mathbb{Q}_\ell$. Recall that $\omega_2:I_\ell\longrightarrow\mathbb{F}_{\ell^2}^\times$ is the
unique character of the inertia subgroup $I_\ell$ given by the rule
$$\tau\longrightarrow\frac{\tau(\ell^{1/(\ell^2-1)})}{\ell^{1/(\ell^2-1)}}.$$
The notation suppresses the dependence on $\ell$, and it would be more appropriate to write
$\omega_{2,\ell}$ instead; the context should be generally clear. One should recall that the
notion of fundamental character is not functorial; the restriction of $\omega_2$ to a local
inertia group $I_v$ is not the second fundamental character of $F_v$ when the ramification degree
of $F_v/\mathbb{Q}_\ell$ is greater than 1. We remark that there is an injection
$\mathbb{F}_{\ell^2}^\times\hookrightarrow\mathrm{GL}_2(\mathbb{F}_\ell)$; it follows that we can
view $\omega_2$ as a 2-dimensional representation $\Omega_2$ over $\mathbb{F}_\ell$. This
representation is irreducible over $\mathbb{F}_\ell$, but if we extend scalars to a coefficient
field of even degree over $\mathbb{F}_\ell$, then $\Omega_2$ becomes reducible, isomorphic over
this quadratic extension to the direct sum of the characters $\omega_2$ and $\omega_2^\ell$.

For an elliptic curve $E$ over a field $F,$ we denote by $E[n]$ the kernel of the multiplication
by $n$ map $E\stackrel{\times n}{\rightarrow}E.$ If $n$ is coprime to the characteristic of $F,$
$$\overline\rho_{E,n}:G_F\longrightarrow
\mathrm{Aut}E[n](\overline{F})\cong\mathrm{GL}_2(\mathbb{Z}/n\mathbb{Z})$$ is the mod~$n$
representation. If $\ell$ is a prime different from the characteristic of $F$, we set
$$\rho_{E,\ell}:G_F\longrightarrow
\lim_{\leftarrow}\mathrm{Aut}E[\ell^n](\overline{F})\cong\mathrm{GL}_2(\mathbb{Z}_\ell).$$

\subsection*{Summary of results}

Let $F$ be a totally real number field, and let $\ell$ be an odd prime. Suppose that for all
$v|\ell$, the ramification index of $F_v/\mathbb{Q}_\ell$ is at most $\ell-1$. Consider
continuous, irreducible representations
$$\rho:\mathrm{Gal}(\overline{F}/F)\longrightarrow
\mathrm{GL}_2(\overline{\mathbb{Q}}_\ell)$$ with determinant the $\ell$-adic cyclotomic
character, and having the same absolutely irreducible residual representation $\overline\rho.$ We
assume that all Artinian quotients of $\rho$ are finite flat at primes above $\ell$, and we
assume further that
$$\overline\rho|_{I_v}\sim
\Omega_2|_{I_v}\quad\mathrm{for~every}\quad v|\ell$$ where $\Omega_2$ is the second fundamental
character of $\mathbb{Q}_\ell$, as in the notation section above, regarded as a 2-dimensional
representation~--~as our coefficient field has residue field containing $\mathbb{F}_{\ell^2}$,
the representation splits as $\omega_2\oplus\omega_2^\ell$. This is the form of the local Galois
representations associated to an elliptic curve with good supersingular reduction at $v$, where
$F_v$ is {\it unramified} over $\mathbb{Q}_\ell$. (If $F_v$ is not unramified, however, the local
Galois representation may take a different form; see section~\ref{modularity} for an example.)
The main applications of the results of the paper will be to such elliptic curves.

Our main result is then:
\begin{theorem} Let $\rho$ be a representation of the above form. Suppose  that
$\overline{\rho}$ has a modular lift which is finite flat at primes above $\ell.$ Assume that
$$\left.\overline\rho\right|_{\mathrm{Gal}\left(\overline{F}/F\left(\zeta_\ell
\right)\right)}$$ is absolutely irreducible, and furthermore assume that \begin{itemize} \item if
$\ell=5$ and $\textrm{Proj}\ \overline\rho|_{\mathrm{Gal}\left(\overline{F}/F\left(\zeta_\ell
\right)\right)}\cong A_5$ then $[F(\zeta_\ell):F]=4.$\end{itemize} Then $\rho$ is also modular.
\end{theorem}

We give two applications of the above. The first relates to Serre's conjecture for mod $7$
representations; we extend the result in~\cite{mano}, and show that:
\begin{theorem}
Let $\overline\rho:\mathrm{Gal}(\overline{\mathbb{Q}}/\mathbb{Q})\longrightarrow
\mathrm{GL}_2(\mathbb{F}_7)$ be an absolutely irreducible, continuous, odd representation.
Suppose that the projective image of inertia at $3$ has odd order and that the determinant of
$\overline\rho$ restricted to the inertia group at $7$ has even order. Then $\overline\rho$ is
modular.
\end{theorem}

This theorem has been used by Dieulefait and the second author~\cite{DM} to give a new criterion
for the modularity of rigid Calabi-Yau threefolds. Of course, it is largely subsumed within
recent work of Khare and Wintenberger; however, we need no hypothesis at~$2$.

Our second application relates to the modularity of elliptic curves over totally real fields. For
general totally real fields, we prove modularity subject to quite a few restrictions. For the
full result, see section~\ref{AppII}. A particularly neat corollary is the following.
\begin{theorem}
Every semistable elliptic curve over $\mathbb{Q}(\sqrt{2})$ is modular.
\end{theorem}
This has implications for the study of certain Diophantine equations, and notably the Fermat
equation, over $\mathbb{Q}(\sqrt{2})$ (see~\cite{JM}).


\section{Local deformations and cohomology groups}\label{conrad}

Our objective in this section is to give good upper bounds on the size of certain local
cohomology groups. We do this for representations of a certain shape (which can be achieved after
an unramified base change). But before that, we begin by setting out our notation. Apart from
$\ell$ being the residue characteristic and $\lambda$ being a uniformizer (instead of $p$ and
$\pi$), our choice of notation is meant to be consistent with \cite{conrad1}.

Throughout this section, we fix a finite field $k$ of characteristic $\ell\geq 3$. We denote by
$A$ its Witt ring $W(k)$ and by $K$ the fraction field of $A$. We fix a finite totally ramified Galois
extension $K'$ of $K$ and denote by $A'$ its ring of integers. We assume that the absolute
ramification index $e=[K':K]$ is less than or equal to $\ell-1.$ The reason for this is that
there is then a good notion of Honda system associated to group schemes.
We also fix throughout a uniformizer $\lambda$ such that $\lambda^e=\epsilon \ell$ with
$\epsilon\in A^\times$ (as $K'$ is a tamely ramified extension). Write $\mathfrak{m}$ for the
maximal ideal of $A'$.

We denote by $\sigma$ the Frobenius automorphism of $A$, and by $D_k$ the Dieudonn\'{e} ring.
Recall that $D_k$ is the $A$-algebra generated by $F$ and $V$ subject to the usual relations
$FV=\ell=VF$, $F\alpha=\sigma (\alpha)F,\ V\alpha=\sigma^{-1}(\alpha)V$ (for $\alpha\in A$). If
there is no cause for confusion, we will abbreviate $D_k$ to simply $D$.

Various tensor products appear in this section. The unspecified $-\otimes -$ will simply mean
$-\otimes_{\mathbb{Z}_\ell}-.$

We shall be working with finite Honda systems over $A'$. For the various properties, see
Conrad~(\cite{conrad1} and~\cite{conrad2}).

We now fix a second finite field $\mathbb{F}$ of characteristic $\ell$ and a continuous
representation
$$\overline{\rho}:G_{K'}\longrightarrow\mathrm{GL}_2(\mathbb{F}).$$
We will shortly impose a further restriction, but for the moment we assume that the
representation is finite---that is, there is a finite flat group scheme over $A'$ whose
associated Galois module (from the generic fibre) gives precisely our representation
$\overline{\rho}.$ This allows us to introduce certain cohomology groups
$H^1_f(G_{K'},\mathrm{ad}\,\overline{\rho})$ and $H^1_f(G_{K'},\mathrm{ad}^0\overline{\rho}).$ We
recall the definitions (see \cite{ddt} for details): elements of
$H^1_f(G_{K'},\mathrm{ad}\,\overline{\rho})$ are the deformations of $\overline{\rho}$ to
$\mathbb{F}[\epsilon]/(\epsilon^2)$ which are finite, and
$H^1_f(G_{K'},\mathrm{ad}^0\overline{\rho})$ is the subspace of
$H^1_f(G_{K'},\mathrm{ad}\,\overline{\rho})$ with determinant (of the deformation) equal to the
determinant of $\overline{\rho}.$

We now impose a restriction on the shape of $\overline{\rho}$:
\begin{assumption}\label{ass1}\rm
$\overline{\rho}$ is equivalent to $\Omega_2|_{G_{K'}}.$\end{assumption}

Let $M$  be the $D_k\otimes\mathbb{F}$-module
$$(k\otimes\mathbb{F})\mathrm{\bf{e}}_1\oplus(k\otimes\mathbb{F})\mathrm{\bf{e}}_2$$
with $F$ and $V$ actions given by
$$\begin{matrix} F(\mathrm{\bf{e}}_1)=0, & F(\mathrm{\bf{e}}_2)=\mathrm{\bf{e}}_1;\\
V(\mathrm{\bf{e}}_1)=0, & V(\mathrm{\bf{e}}_2)=-\mathrm{\bf{e}}_1.\end{matrix}$$ (To be more
precise, these give the action on our basis elements which one then extends Frobenius
semi-linearly.) Let $L$ be the subspace $(k\otimes\mathbb{F})\mathrm{\bf{e}}_2.$ Then $(L,M)$ is
the finite Honda system over $A$ associated to $\Omega_2|_{G_K}.$ This follows, after base change
(see Section 4 of \cite{conrad1}), from the description of the Honda system over
$\mathbb{Z}_\ell$ associated to $\Omega_2$. (This is presumably well known, but a proof is given
in the appendix.) We reserve $(L,M)$ for this particular Honda system throughout.

By the results of \cite{conrad1}, calculating $H^1_f(G_{K'},\mathrm{ad}\,\overline{\rho})$ is the
same as calculating extensions of $(L,M)$ by itself in the category of finite Honda systems over
$A'.$ As a first step to this calculation, we investigate the extensions of $M$ by itself in the
category of $D_k\otimes\mathbb{F}$ modules.

We begin with a technical lemma which enables us to reduce calculations to one of linear algebra.
\begin{lemma}\label{freeness} Let $R$ be a ring with finite
cardinality. If $$0\longrightarrow R^m\longrightarrow U\longrightarrow R^n\longrightarrow 0$$ is
an exact sequence of $R-$modules, then $U$ is free and isomorphic to $R^{n+m}.$\end{lemma}

\begin{Proof} The exact sequence implies that $U$ can be generated
by $n+m$ elements. Hence there is a surjective $R$-module homomorphism $R^{n+m}\twoheadrightarrow
U.$ As $R$ has finite cardinality, we get $R^{n+m}\cong U.$\end{Proof}

\begin{proposition} \label{ext1} The group of extensions
$\mathrm{Ext}^1_{D_k\otimes\mathbb{F}}(M,M)$ is (non-canonically) isomorphic as an
$\mathbb{F}$-vector space to
\begin{itemize}
\item $(k\otimes\mathbb{F})\oplus (\mathbb{F}_\ell\otimes\mathbb{F})$ if the degree
$[k:\mathbb{F}_\ell]$ is odd, and \item $(k\otimes\mathbb{F})\oplus
(\mathbb{F}_{\ell^2}\otimes\mathbb{F})$ if the degree $[k:\mathbb{F}_\ell]$ is even.
\end{itemize}
\end{proposition}

\begin{Proof}
By Lemma~\ref{freeness}, we can certainly take any extension class, as an $A\otimes\mathbb{F}$
module, to be
$$M\oplus M=
\Big( (k\otimes\mathbb{F})(\mathrm{\bf{e}}_1,0)\oplus(k\otimes\mathbb{F})(\mathrm{\bf{e}}_2,0)
\Big)\oplus \Big(
(k\otimes\mathbb{F})(0,\mathrm{\bf{e}}_1)\oplus(k\otimes\mathbb{F})(0,\mathrm{\bf{e}}_2) \Big).$$
We need to specify the actions of $F$ and $V$. In order to do this, we write down matrices using
the above choice of basis and compute (remembering to keep track of Frobenius semi-linearity).

To begin with, we can write
$$F=\begin{pmatrix} 0 & 1 & f_1 & f_2\\
0& 0 & f_3 & f_4\\
0 & 0 & 0 & 1\\
0 & 0 & 0 & 0\end{pmatrix}\ \mathrm{and}\
V=\begin{pmatrix} 0 & -1 & v_1 & v_2\\
0& 0 & v_3 & v_4\\
0 & 0 & 0 & -1\\
0 & 0 & 0 & 0\end{pmatrix}.$$ Since $FV=VF=\ell=0,$ we must have the following equalities:
\begin{eqnarray*}
\begin{pmatrix} 0&1\\ 0&0\end{pmatrix}
\begin{pmatrix}\sigma (v_1) &\sigma (v_2)\\ \sigma (v_3) &\sigma
(v_4)\end{pmatrix} +
\begin{pmatrix} f_1 & f_2\\ f_3 & f_4\end{pmatrix}
\begin{pmatrix} 0&-1\\ 0&0\end{pmatrix}
& = & 0\\
\begin{pmatrix} 0&-1\\ 0&0\end{pmatrix}
\begin{pmatrix}\sigma^{-1} (f_1) &\sigma^{-1} (f_2)\\ \sigma^{-1} (f_3)
&\sigma^{-1} (f_4)\end{pmatrix} +
\begin{pmatrix} v_1 & v_2\\ v_3 & v_4\end{pmatrix}
\begin{pmatrix} 0&1\\ 0&0\end{pmatrix} & = & 0
\end{eqnarray*}
Multiplying out, we find that
$$f_3=v_3=0,\ \mathrm{and}\ f_1=\sigma(v_4), f_4=\sigma(v_1).$$

We now reduce the number of variables further by applying appropriate $k\otimes\mathbb{F}$-linear
automorphisms of $M\oplus M.$ Let $A$ be the endomorphism
$$\begin{pmatrix}
1&0&a_1&a_2\\0&1&a_3&a_4\\0&0&1&0\\0&0&0&1\end{pmatrix}.$$ To calculate $AFA^{-1},$ we need to
calculate
$$\begin{pmatrix}a_1&a_2\\a_3&a_4\end{pmatrix}\begin{pmatrix}0&1\\0&0\end{pmatrix}
-\begin{pmatrix}0&1\\0&0\end{pmatrix}
\begin{pmatrix}\sigma(a_1)&\sigma(a_2)\\\sigma(a_3)&\sigma(a_4)\end{pmatrix}
+\begin{pmatrix}f_1&f_2\\0&f_4\end{pmatrix}$$ which is
$$\begin{pmatrix} -\sigma(a_3)&
a_1-\sigma(a_4)\\0&a_3\end{pmatrix} +\begin{pmatrix}f_1&f_2\\0&f_4\end{pmatrix}.$$ We can thus
assume that $f_4=f_2=0,$ which implies that $v_1=0.$ Under this assumption, our choice of $A$ is
then restricted to
$$a_3=0 \ \mathrm{and}\ a_1=\sigma(a_4).$$

To calculate $AVA^{-1},$ we need to compute
$$\begin{pmatrix}a_1&a_2\\0&a_4\end{pmatrix}\begin{pmatrix}0&-1\\0&0\end{pmatrix}
-\begin{pmatrix}0&-1\\0&0\end{pmatrix}
\begin{pmatrix}\sigma^{-1}(a_1)&\sigma^{-1}(a_2)\\0&\sigma^{-1}(a_4)\end{pmatrix}
+\begin{pmatrix}0&v_2\\0&v_4\end{pmatrix}$$ which is
$$\begin{pmatrix} 0&
-a_1+\sigma^{-1}(a_4)\\0&0\end{pmatrix} +\begin{pmatrix}0&v_2\\0&v_4\end{pmatrix}.$$ Since we
have $a_1=\sigma(a_4),$ our choice $v_2\in k\otimes\mathbb{F}$ can further be restricted to a
choice of representative of an element of
$$\frac{k\otimes\mathbb{F}}{(\sigma^2-1)(k\otimes\mathbb{F})},$$
while $v_4$ can be chosen to be an arbitrary element of $k\otimes\mathbb{F}$. The proposition
then follows.
\end{Proof}

\begin{theorem} The dimension of
$H^1_f(G_{K'},\mathrm{ad}\,\overline{\rho})$ as an $\mathbb{F}$-vector space is at most
\begin{itemize}
\item $[K':\mathbb{Q}_\ell]+2$ if $[k:\mathbb{F}_\ell]$ is even, and \item
$[K':\mathbb{Q}_\ell]+1$ if $[k:\mathbb{F}_\ell]$ is odd.
\end{itemize}
\end{theorem}

\begin{Proof} As in \cite{conrad2}, we have an $\mathbb{F}$-linear map of
vector spaces
$$t:H^1_f(G_{K'},\mathrm{ad}\,\overline{\rho})\longrightarrow\mathrm{Ext}^1(M,M).$$
In words, the map $t$ is just `take Dieudonn\'{e} module of the special fibre of the associated
finite flat group scheme'. We already have a bound for the Ext-group, thanks to
Proposition~\ref{ext1}. We now start analysing the kernel of the above linear map.

We begin by describing the structure of the $A'$-module $M_{A'}.$ We recall the definition (due
to Fontaine), and refer to \cite{conrad1} for the explicit description we need (see Definition
2.1 of \cite{conrad1}). As already set out in the beginning of this section, we have a fixed
uniformizer $\lambda$ of $A'$ satisfying $\lambda^e=\epsilon \ell$ with $\epsilon\in A^\times.$

We have the standard identification of $M^{(1)}=(A,\sigma)\otimes_AM$ with $M$ as an abelian
group and twisted $A$-action. The Dieudonn\'{e} module structure then gives us two $A$-linear
maps
$$F_0:M^{(1)}\longrightarrow M\quad \mathrm{and} \quad V_0:M\longrightarrow
M^{(1)}.$$ (As in \cite{conrad1}, we shall not abbreviate these to $F$ and $V.$) There are
$A'$-linear maps
$$F^{M}:A'\otimes_A M^{(1)}\longrightarrow A'\otimes_A M\quad
\mathrm{and} \quad V^M:\mathfrak{m}\otimes_A M\longrightarrow
\ell^{-1}\mathfrak{m}\otimes_AM^{(1)}$$ obtained simply by tensoring with the identity map on
$A'$ and the map $x\rightarrow \ell^{-1}x$ respectively.

The $A'$-module $M_{A'}$ is then the quotient of
$$(A'\otimes_AM)\oplus(\ell^{-1}\mathfrak{m}\otimes_AM^{(1)})$$
by the submodule
$$\left\{\left(
\phi_0^M(u)-F^M(w),\phi_1^M(w)-V^M(u)\right)\,\left. \right|\,u\in\mathfrak{m}\otimes_AM,w\in
A'\otimes_AM^{(1)}\right\}$$ where $\phi_0^M,\ \phi_1^M$ are the maps $$
\phi_0^M:\mathfrak{m}\otimes_AM\longrightarrow A'\otimes_AM \quad\mathrm{and}\quad
\phi_1^M:A'\otimes_AM^{(1)}\longrightarrow \ell^{-1}\mathfrak{m}\otimes_AM^{(1)}$$ induced by the
inclusions $\mathfrak{m}\hookrightarrow A'$ and $A'\hookrightarrow \ell^{-1}\mathfrak{m}.$

A basis of $A'\otimes_AM$ as a free $k\otimes\mathbb{F}$-module is given by
$$\lambda^i\otimes\mathrm{\bf{e}}_j, \ i=0,\ldots e-1, j=1,2.$$ For
$\ell^{-1}\mathfrak{m}\otimes_AM^{(1)},$ we have the $k\otimes\mathbb{F}$ basis
$$\lambda^{-i}\otimes\mathrm{\bf{e}}_j,\ i=0,1,\ldots e-1, j=1,2.$$

Note that for $i\geq 1,$ the elements $(\lambda^i\otimes\mathrm{\bf{e}}_1,0)$ are trivial in
$M_{A'}.$ Indeed, we have
$$(\lambda^i\otimes\mathrm{\bf{e}}_1,0)=\left(
\phi_0^M(\lambda^i\otimes\mathrm{\bf{e}}_1)-F^M(0),
0-V^M(\lambda^i\otimes\mathrm{\bf{e}}_1)\right).$$ Furthermore, for $i\geq 1,$ we have
\begin{eqnarray*}
(0,\lambda^{-i}\otimes\mathrm{\bf{e}}_1)&=&(0,0-V^M(\lambda^{e-i}\otimes\mathrm{\bf{e}}_2))\\
&=&(-\lambda^{e-i}\otimes\mathrm{\bf{e}}_2,0).
\end{eqnarray*}
Note also that
\begin{eqnarray*}
(0,1\otimes\mathrm{\bf e}_1)&=&\left(\phi_0^M(0)-F^M(1\otimes\mathrm{\bf
e}_1), \phi_1^M(1\otimes\mathrm{\bf e}_1)-V^M(0)\right),\quad\mathrm{and}\\
(0,1\otimes\mathrm{\bf e}_2)&=& (1\otimes\mathrm{\bf e}_1,0)+
\left(\phi_0^M(0)-F^M(1\otimes\mathrm{\bf e}_2), \phi_1^M(1\otimes\mathrm{\bf
e}_2)-V^M(0)\right).
\end{eqnarray*}
Thus any element in $M_{A'}$ can be expressed as an $k\otimes\mathbb{F}$-linear combination of
$$(1\otimes\mathrm{\bf{e}}_1,0),
(\lambda^i\otimes\mathrm{\bf{e}}_2,0)\ \mathrm{and}\ (0,\lambda^{-m}\otimes\mathrm{\bf{e}}_2)$$
with $i=0,1,\ldots e-1$ and $m=1,\ldots e-1.$ Since the $A'$-length of $M_{A'}$ is the same as
the $A$-length of $M$ times $e$ (Lemma 2.2 of \cite{conrad1}), we deduce that the set of
generators above is in fact a basis.

Obviously, the $A'$-submodule of $M_{A'}$ obtained by taking the $A'$-span of $L$ is precisely
$A'\otimes_AA\otimes\mathbb{F}(\mathrm{\bf e}_2,0).$ Now let $(L',M')$ be the finite Honda system
for an element in the kernel of $t.$ Since $M'=M\oplus M$ as a $D_k\otimes \mathbb{F}$-module, we
can write $M'_{A'}=M_{A'}\oplus M_{A'}.$ We must therefore have, by length considerations,
$$L'=(A'\otimes_AA\otimes\mathbb{F})((\mathrm{\bf e}_2,0),0)+
(A'\otimes_AA\otimes\mathbb{F}) (x,(\mathrm{\bf e}_2,0))$$ for some $x\in M_{A'}.$ From our
description of a basis of $M_{A'},$ it follows that we can take
$$x=a(1\otimes\mathrm{\bf{e}}_1,0)+y$$
with $a\in k\otimes\mathbb{F}$ and $y$ an element in the $A\otimes\mathbb{F}$-span of
$(0,\lambda^{-m}\otimes\mathrm{\bf{e}}_2),~m=1,\ldots e-1.$ By applying a
$D_k\otimes\mathbb{F}$-linear automorphism of $M\oplus M$ of the type
$$\begin{pmatrix}
1&0&0&*\\0&1&0&0\\0&0&1&0\\0&0&0&1\end{pmatrix},$$ we can assume that $a=0.$ Hence the kernel has
dimension, as an $\mathbb{F}$-vector space, at most $(e-1)[k:\mathbb{F}_\ell]$; and this proves
the theorem.
\end{Proof}

\begin{corollary}\label{cor2.5}
The dimension of $H^1_f(G_{K'},\mathrm{ad}^0\overline{\rho})$ as an $\mathbb{F}$-vector space is
at most
\begin{itemize}
\item $[K':\mathbb{Q}_\ell]+1$ if $[k:\mathbb{F}_\ell]$ is even, and \item $[K':\mathbb{Q}_\ell]$
if $[k:\mathbb{F}_\ell]$ is odd.
\end{itemize}
\end{corollary}

\section{The deformation problem}\label{sec3}

We now set up the deformation problem we want to study. We begin by fixing a totally real
extension $F$ of even degree (over $\mathbb{Q}$), an odd prime $\ell$, a finite field $k$ of
characteristic $\ell$, and a continuous homomorphism
$$\overline{\rho}:G_F\longrightarrow\mathrm{GL}_2(k)$$ which is absolutely
irreducible and odd. We assume that the ramification degree of $F$ at all primes over $\ell$ is
less than or equal to $\ell-1.$ Further, we suppose that $\overline{\rho}$ has the following
properties:
\begin{itemize}
\item The determinant of $\overline{\rho}$ is the mod $\ell$ cyclotomic character. \item
$\overline{\rho}$ restricted to the absolute Galois group of $F(\zeta_\ell)$ is absolutely
irreducible. \item If $\ell=5$ and $\textrm{Proj}\
\overline\rho|_{\mathrm{Gal}\left(\overline{F}/F\left(\zeta_\ell \right)\right)}$ then
$[F(\zeta_\ell):F]=4.$ \item Let $x$ be a prime of $F$ above $\ell$ and let $I_x$ the inertia
group of $F_x.$ Then
$$\overline\rho |_{I_x}\sim\Omega_2|_{I_x}$$
where $\Omega_2$ is the second fundamental character.
\end{itemize}
We assume that the characteristic polynomial of $\overline\rho (\sigma)$ is split over $k$ for
any $\sigma\in G_F$. We fix a finite extension $K$ of $\mathbb{Q}_\ell$ with ring of integers
$\mathcal{O},$ maximal ideal $(\lambda)$ and residue field $k.$

Let $\mathcal{C}_\mathcal{O}$ be the category of complete, local, Noetherian
$\mathcal{O}$-algebras with residue field $k.$ Given
$(A,\mathfrak{m}_A)\in\mathcal{C}_\mathcal{O},$ we call a continuous homomorphism
$$\rho_A:G_F\longrightarrow\mathrm{GL}_2(A)$$ a {\it finite flat deformation} of
$\overline\rho$ if
\begin{itemize}
\item $\rho_A$ is odd and unramified outside finitely many primes, \item
$\rho_A\pmod{\mathfrak{m}_A}=\overline\rho$, \item $\rho_A$ is finite flat at primes $v|\ell$
(i.e., the restriction of $\rho_A$ to $G_{F_v}$, for $v|\ell$, has the property that for all
$n\ge1$, the $F_v$-group scheme associated to the $G_{F_v}$-module $\rho_A\mbox{ mod
}\mathfrak{m}_A^n$ is the generic fibre of a finite flat group scheme over $\mathcal{O}_{F,v}$),
and \item $\rho_A$ has determinant the $\ell$-adic cyclotomic character.
\end{itemize}
Two such deformations are said to be strictly equivalent if one can be conjugated to the other by
a matrix which reduces to the identity modulo the maximal ideal $\mathfrak{m}_A.$

Now let $\Sigma$ be a finite set of (finite) primes of $F$ not containing any places over $\ell$
(and it could be empty). We say a finite flat deformation is of type $\Sigma$ if the
representation is unramified outside primes in $\Sigma$ and outside the set of primes where
$\overline\rho$ is ramified. There is then a universal finite flat deformation of $\overline\rho$
of type $\Sigma$ which we shall denote by $(R_\Sigma,\rho_\Sigma).$

Given a finite flat deformation $\rho:G_F\rightarrow \mathrm{GL}_2(\mathcal{O}/\lambda^n)$ of
type $\Sigma ,$ one defines the Galois cohomology group $H^1_\Sigma(G_F,\mathrm{ad}^0\rho)$ to be
the deformations of $\rho$ to $(\mathcal{O}/\lambda^n)[\epsilon]/\epsilon^2$ which are of type
$\Sigma.$ Recall that $\mathrm{ad}^0\rho$ can be identified with the group of $2\times2$ trace
zero matrices over $\mathcal{O}/\lambda^n$ with $G_F$ action via conjugation (by $\rho$). The
cohomology group $H^1_\Sigma (G_F, \mathrm{ad}^0\rho)$ is then precisely
$H^1_{\mathcal{L}_\Sigma} (G_F, \mathrm{ad}^0\rho)$ where the local conditions
$\mathcal{L}_\Sigma=\{L_x\}$ are given by:
\begin{itemize}
\item $L_x=H^1(G_{F_x}/I_x, \mathrm{ad}^0\rho^{I_x})$ if $x\nmid\ell,$ $x\notin\Sigma$ and
$\overline\rho$ is unramified at $x,$ \item $L_x=H^1(G_{F_x},\mathrm{ad}^0\rho)$ if $x\nmid\ell
,$ and either $x\in\Sigma$ or $\overline\rho$
 is ramified at $x,$
\item $L_x=H^1_f(G_{F_x},\mathrm{ad}^0\rho)$ if $x|\ell.$
\end{itemize}

The universal deformation ring $R_\Sigma$ can be topologically generated as an
$\mathcal{O}$-algebra by $\mathrm{dim}_kH^1_\Sigma(G_F,\mathrm{ad}^0\overline\rho)$ elements. If
$\pi:R_\Sigma\twoheadrightarrow\mathcal{O}$ is an $\mathcal{O}$-algebra homomorphism with
corresponding representation $\rho,$ we have a canonical isomorphism
$$\mathrm{Hom}\left(\mathrm{ker}\,\pi/(\mathrm{ker}\,\pi)^2,K/\mathcal{O}\right)\cong
H^1_\Sigma(G_F,\mathrm{ad}^0\rho\otimes K/\mathcal{O}).$$

The pairing $\mathrm{ad}^0\overline\rho\times\mathrm{ad}^0\overline\rho\rightarrow k$ obtained by
taking the trace is perfect. Using this pairing, one defines
$H^1_\Sigma(G_F,\mathrm{ad}^0\overline\rho(1))$ to be given by local conditions $\{L_x^\bot\}$
where $L_x^\bot$ is the orthogonal complement to $L_x$ with respect to the perfect pairing
$$H^1(G_{F_x},\mathrm{ad}^0\overline\rho)\times H^1(G_{F_x},\mathrm{ad}^0\overline\rho(1))
\longrightarrow H^2(G_{F_x},k(1))\simeq k.$$

From now onwards, we assume the following:
\begin{assumption}\label{ass2}\rm
For each prime $x$ of $F$ dividing $\ell$, the Honda system associated to $\overline\rho|_{F_x}$
has the particular form specified in Assumption~\ref{ass1}.
\end{assumption}

Now we make some calculations of these cohomology groups, using similar arguments to those of
Wiles.
\begin{theorem}\label{th3.2-1}
As an $\mathcal{O}$-algebra,
$$\mathrm{dim}_kH^1_\Sigma(G_F,\mathrm{ad}^0\overline\rho(1))+\sum_{x\in\Sigma}\mathrm{dim}_k
H^0(G_{F_x},\mathrm{ad}^0\overline\rho(1))$$ elements are sufficient to generate the universal
deformation ring $R_\Sigma$ topologically.
\end{theorem}

\begin{Proof} This is almost exactly the same as the proof of Corollary 2.43 in
\cite{ddt}. Using Theorem 2.19 of \cite{ddt} (a full proof is given in \cite{nsw}, p.440), one
finds that $\mathrm{dim}_kH^1_\Sigma(G_F,\mathrm{ad}^0\overline\rho)$ is the sum of terms:
\begin{itemize}
\item $\mathrm{dim}_kH^1_\Sigma(G_F,\mathrm{ad}^0\overline\rho(1))$; \item
$\sum_{x|\ell}\mathrm{dim}_kH^1_f(G_{F_x})- \sum_{x|\ell}\mathrm{dim}_kH^0(G_{F_x})-
\sum_{x|\infty}\mathrm{dim}_kH^0_\Sigma(G_{F_x}),$ where $H^*_*(G_{F_x})$ means the cohomology
group $H^*_*(G_{F_x},\mathrm{ad}^0\overline\rho)$. This term is less than or equal to $0$ by
Corollary~\ref{cor2.5}. \item $\mathrm{dim}_kH^1(G_{F_x},\mathrm{ad}^0\overline\rho)-
\mathrm{dim}_kH^0(G_{F_x},\mathrm{ad}^0\overline\rho),$ which equals
$\mathrm{dim}_kH^0(G_{F_x},\mathrm{ad}^0\overline\rho(1)),$ for each $x\in\Sigma.$
\end{itemize}
\end{Proof}

Theorem 2.49 of \cite{ddt} still holds in our present setting; the proof, with trivial
modifications, remains valid. The result being of significant importance, we give a brief sketch
of the proof.

\begin{theorem}\label{th3.2-2}
Let $r= \mathrm{dim}_kH^1_\emptyset(G_F,\mathrm{ad}^0\overline\rho(1))$. For every positive
integer $n,$ we can find a finite set primes $\Sigma_n$ such that the following hold:
\begin{itemize}
\item Every prime in $\Sigma_n$ has norm congruent to $1$ modulo $\ell^n$; \item The sets
$\Sigma_n$ all have size equal to $r$; \item If $x\in\Sigma_n,$ then $\overline\rho$ is
unramified at $x$ and the Frobenius (at $x$) has distinct eigenvalues; \item The universal
deformation ring $R_{\Sigma_n}$ can be topologically generated as an $\mathcal{O}$-algebra by $r$
elements.
\end{itemize}
\end{theorem}

\begin{Proof} As in the proof of Theorem 2.49 of
\cite{ddt}, one reduces the result to showing that for $\psi\in H^1_\emptyset
(G_F,\mathrm{ad}^0\overline\rho (1))-\{0\}$, we can find a $\sigma\in G_F$ such that
\begin{itemize}
\item $\sigma$ acts trivially on $F(\zeta_{\ell^n})$, \item $\mathrm{ad}^0\overline\rho(\sigma)$
has an eigenvalue not equal to 1, and \item
$\psi(\sigma)\notin(\sigma-1)\mathrm{ad}^0\overline\rho(1).$
\end{itemize}
(We remark that Theorem~\ref{th3.2-1} is crucial in getting the right number of generators from
this reduction.)

Let $F_n$ be the minimal extension of $F(\zeta_{\ell^n})$ on which $\mathrm{ad}^0\overline\rho$
acts trivially. The degree of the extension $F_1/F_0$ is at most $\ell-1$; the degree $[F_n:F_1]$
is of $\ell$-power order. It follows that
$$H^1(\mathrm{Gal}(F_n/F_0),\mathrm{ad}^0\overline\rho(1))^{G_F}\cong
\mathrm{Hom}(\mathrm{Gal}(F_n/F_1),\mathrm{ad}^0\overline\rho(1)^{G_F})$$ is trivial (since
$\overline{\rho}$ restricted to the absolute Galois group of $F(\zeta_\ell)$ is absolutely
irreducible).

Now consider $H^1(\mathrm{Gal}(F_0/F), \mathrm{ad}^0\overline\rho(1)^{G_{F_0}})$. If this is
non-trivial, the order of $\mathrm{Gal}(F_0/F)$ must be divisible by $\ell$ and
$\mathrm{Gal}(F_0/F)$ must have $\mathrm{Gal}(F(\zeta_\ell)/F)$ as a quotient. Note that
$\mathrm{Gal}(F_0/F)$ is isomorphic to the projective image of $\overline\rho,$ and so from the
list in Theorem 2.47 of \cite{ddt} we see that the case $\ell=5$ and $\textrm{Proj}\
\overline\rho|_{\mathrm{Gal}\left(\overline{F}/F\left(\zeta_\ell \right)\right)}$ cannot occur.
In the other cases the projective image of $\overline\rho$ is a semi-direct extension of
$PSL_2(\mathbb{F}_{{\ell}^r})$ by a group of order prime to $\ell,$ and so
$H^1(\mathrm{Gal}(F_0/F),\mathrm{ad}^0\overline\rho(1))$ again vanishes on applying Lemma 2.48 of
\cite{ddt}.

A straightforward application of the inflation-restriction sequence then implies that the group
$H^1(\mathrm{Gal}(F_n/F), \mathrm{ad}^0\overline\rho(1))$ is trivial, and it follows that
$\psi(G_{F_n})$ is non-trivial.

Now $\overline\rho$ restricted to $G_{F(\zeta_{\ell^n})}$ is still absolutely irreducible. Thus
the order of $\mathrm{Gal}(F_n/F(\zeta_{\ell^n}))$ is not a power of $\ell$. The group
$\mathrm{Gal}(F_n/F(\zeta_{\ell^n}))$ also acts (non-trivially) on
$\{0\}\neq\psi(G_{F_n})\subset\mathrm{ad}^0\overline\rho$. Therefore we can find a non-trivial
element $g\in\mathrm{Gal}(F_n/F(\zeta_{\ell^n}))$ of order prime to $\ell$ and fixing a non-zero
element of $\psi(G_{F_n})$. Let $\tilde{g}\in G_{F(\zeta_{\ell^n})}$ be a lift of $g$. As
$\psi(G_{F_n})\not\subset(g-1)\mathrm{ad}^0\overline\rho(1)$, we can find an $h\in G_{F_n}$ such
that
$$\psi (h\tilde{g})=\psi (h)+\psi
(\tilde{g})\notin(\tilde{g}-1)\mathrm{ad}^0\overline\rho(1).$$

Finally, take $\sigma=h\tilde{g}$. Then $\sigma$ acts trivially on $F(\zeta_{\ell^n}),$ and
$(\sigma-1)\mathrm{ad}^0\overline\rho(1)=
\tilde{g}-1)\mathrm{ad}^0\overline\rho(1)\not\supset\psi(\sigma)$. Since the order of $\sigma$ is
prime to $\ell$ (and is not 1), it follows that $\mathrm{ad}^0\overline\rho(\sigma)$ has an
eigenvalue not equal to~1.
\end{Proof}

\section{Hecke algebras and $\ell$-adic modular forms}
We fix a totally real field $F$ of even degree and an odd rational prime $\ell$. We write $D$ for
the division algebra with centre $F$ and ramified exactly at the set of infinite places of $F$.
Write $Z$ for the algebraic group defined by $Z(R)=(D\otimes_FR)^\times$ if $R$ is an
$F$-algebra. We also fix the following:
\begin{itemize}
\item A maximal order $\mathcal{O}_D,$ and isomorphisms $\mathcal{O}_{D,x}\cong
M_2(\mathcal{O}_{F,x})$ for all finite places $x$ of $F.$ These isomorphisms give us an
identification of $\mathrm{GL}_2(\mathbb{A}_F^\infty)$ with
$(D\otimes_\mathbb{Q}\mathbb{A}^\infty)^\times$; \item A uniformiser $\varpi_x$ of
$\mathcal{O}_{F,x}$ for each finite place $x$.
\end{itemize}
We write $A$ for a topological $\mathbb{Z}_\ell$-algebra which is one of the following: a finite
extension of $\mathbb{Q}_\ell,$ the ring of integers in such an extension, or a quotient of such
a ring of integers.

\begin{definition}\rm
For a compact open subgroup $U\subset(D\otimes_\mathbb{Q}\mathbb{A}^\infty)^\times$ and a
topological ring $A$ as above, we define $S_A(U)$ to be the space of continuous functions
$$f:D^\times\backslash(D\otimes_\mathbb{Q}\mathbb{A}^\infty)^\times/U.Z(\mathbb{A}_F^\infty)\longrightarrow A.$$
We define $S_A$ to be the direct limit of $S_A(U)$ as $U$ varies over open compact subsets of
$(D\otimes_\mathbb{Q}\mathbb{A}^\infty)^\times$.
\end{definition}

For a compact open $U$, the finite double coset decomposition
$$(D\otimes_\mathbb{Q}\mathbb{A}^\infty)^\times=\coprod D^\times
t_iU.Z(\mathbb{A}_F^\infty)$$ shows that
\begin{eqnarray*}
S_A(U)&\longrightarrow&\bigoplus_iA\\
f&\longrightarrow&(f(t_i))_i
\end{eqnarray*}
is an isomorphism. In particular, for any $A$-algebra $B,$ we have
$$S_A(U)\otimes_AB\cong S_B(U).$$
We denote by $[t_i]$ the function in $S_A(U)$ which is $1$ on $D^\times
t_iU.Z(\mathbb{A}_F^\infty)$ and $0$ elsewhere.

\begin{definition}\rm
For an ideal $\mathfrak{n}$ of $\mathcal{O}_F$ and quotients $H_x$ of
$(\mathcal{O}_{F,x}/\mathfrak{n}_x)^\times,$ we set $H=\prod_xH_x$. We define $U_H(\mathfrak{n})$
to be the compact open subgroup $\prod_xU_H(\mathfrak{n})_x\subset
(D\otimes_\mathbb{Q}\mathbb{A}^\infty)^\times$ where
$$U_H(\mathfrak{n})_x=\left.\left\{
\begin{pmatrix}a&b\\c&d\end{pmatrix}\in
\mathrm{GL}_2(\mathcal{O}_{F,x})\cong\mathcal{O}_{D,x}^\times\,\right|\,c\in\mathfrak{n}_x,\
ad^{-1}=1\ \mathrm{in}\ H_x\right\}.$$
\end{definition}

Now let $\mathfrak{n}$ and $H_x$ be as in the above definition. We recall the definitions of the
various Hecke operators on $S_A(U_H(\mathfrak{n}))$:
\begin{itemize}
\item If $x$ does not divide $\ell\mathfrak{n},$ we denote the Hecke operators
$$\left[U_H(\mathfrak{n})\begin{pmatrix}
\varpi_x&0\\0&1\end{pmatrix}U_H(\mathfrak{n})\right]\quad \mathrm{and}\quad
\left[U_H(\mathfrak{n})\begin{pmatrix} \varpi_x&0\\0&\varpi_x
\end{pmatrix}U_H(\mathfrak{n})\right]$$ by $T_x$ and $S_x$ respectively. \item If $x$ divides
$\mathfrak{n},$ we set
$$\langle h\rangle=\left[U_H(\mathfrak{n})\begin{pmatrix}
\tilde{h}&0\\0&1\end{pmatrix}U_H(\mathfrak{n})\right]$$ for $h\in H_x$ and $\tilde{h}$ a choice
of lift of $h$ to $\mathcal{O}_{F,x}^\times$. \item If $x$ divides $\mathfrak{n},$ the Hecke
operators
$$ \left[U_H(\mathfrak{n})\begin{pmatrix}
\varpi_x&0\\0&1\end{pmatrix}U_H(\mathfrak{n})\right]\quad\mathrm{and}\quad
\left[U_H(\mathfrak{n})\begin{pmatrix} 1&0\\0&\varpi_x\end{pmatrix}U_H(\mathfrak{n})\right]$$ are
denoted by $\mathrm{\bf U}_{\varpi_x}$ and $\mathrm{\bf V}_{\varpi_x}$ respectively. We also
denote by $S_x$ the Hecke operator
$$\left[U_H(\mathfrak{n})\begin{pmatrix}
\varpi_x&0\\0&\varpi_x \end{pmatrix}U_H(\mathfrak{n})\right].$$
\end{itemize}

\begin{definition}\rm
Let $\mathfrak{n},\ H_x$ and $A$ be as in the preceding paragraphs. We define the {\sl Hecke
algebra} $\mathbb{T}_A(U_H(\mathfrak{n}))$ to be the $A$-subalgebra of
$\mathrm{End}_A(S_A(U_H(\mathfrak{n})))$ generated by $T_x$ (for $x$ not dividing
$\ell\mathfrak{n}$) and $\mathrm{\bf U}_{\varpi_x}$ (for $x|\mathfrak{n}$ but not dividing
$\ell$).

A maximal ideal $\mathfrak{m}$ of $\mathbb{T}_A(U_H(\mathfrak{n}))$ is said to be {\sl
Eisenstein} if it contains $T_x-2$ and $S_x-1$ for all but finitely many primes with $\mathrm{\bf
N}x\pmod{\ell}=1$.
\end{definition}

The Hecke algebra $\mathbb{T}_A(U_H(\mathfrak{n}))$ is always commutative. Also,
$\mathbb{T}_{\mathbb{Z}_\ell}(U_H(\mathfrak{n}))$ is semi-local and $\ell$-adically complete, and
we have the identification
$$\mathbb{T}_{\mathbb{Z}_\ell}(U_H(\mathfrak{n}))\cong
\prod\mathbb{T}_{\mathbb{Z}_\ell}(U_H(\mathfrak{n}))_\mathfrak{m}$$ where the product is over all
maximal ideals $\mathfrak{m}.$

If either $\ell$ is invertible in $A,$ or if $\mathbb{Q}(\zeta +\zeta^{-1})\not\subset F$ where
$\zeta$ is a primitive $\ell$th root of unity, we have a perfect pairing on
$S_A(U_H(\mathfrak{n}))$ defined by
$$(f_1,f_2)_{U_H(\mathfrak{n})}=\sum_i
f_1(t_i)f_2(t_i)\left(\#\frac{U_H(\mathfrak{n}).Z(\mathbb{A}_F^\infty) \cap t_i^{-1}D^\times
t_i}{F^\times}\right)^{-1}$$ where
$$(D\otimes_\mathbb{Q}\mathbb{A}^\infty)^\times=\coprod D^\times
t_iU_H(\mathfrak{n}).Z(\mathbb{A}_F^\infty).$$ We call this the {\sl standard pairing}. The Hecke
operators are not necessarily self-adjoint with respect to this pairing; the general behaviour of
operators is given by
$$\left(
\left[U_{H'}(\mathfrak{n}')gU_H(\mathfrak{n})\right]f_1,f_2\right)_{U_{H'}(\mathfrak{n}')}=
\left(f_1,
\left[U_H(\mathfrak{n})g^{-1}{U_{H'}}(\mathfrak{n}')\right]f_2\right)_{U_H(\mathfrak{n})}.$$

Now fix a finite set of primes $\Sigma$, none lying above $\ell$, and let
$\mathfrak{n}_\Sigma=\prod_{x\in\Sigma}x^2$. Let $K$ be a finite extension of $\mathbb{Q}_\ell$
which contains all embeddings $F\hookrightarrow\overline{\mathbb{Q}}_\ell$, and let $\mathcal{O}$
be its ring of integers. We fix a decomposition
$$(D\otimes_\mathbb{Q}\mathbb{A}^\infty)^\times=\coprod D^\times
g_iU_{1}(\mathfrak{n}_\Sigma).Z(\mathbb{A}_F^\infty)\amalg\coprod D^\times
h_iU_{1}(\mathfrak{n}_\Sigma).Z(\mathbb{A}_F^\infty)$$ where the $g_i$'s and $h_i$'s are such
that
$$\ell\not\left|\#\frac{U_1(\mathfrak{n}_\Sigma).Z(\mathbb{A}_F^\infty)
\cap g_i^{-1}D^\times g_i}{F^\times}\right.\quad\mathrm{and}\quad
\ell\left|\#\frac{U_1(\mathfrak{n}_\Sigma).Z(\mathbb{A}_F^\infty) \cap h_i^{-1}D^\times
h_i}{F^\times}\right..$$ We denote by $S_\mathcal{O}(U_1(\mathfrak{n}_\Sigma))^*$ the
$\mathcal{O}$-submodule of $S_\mathcal{O}(U_1(\mathfrak{n}_\Sigma))$ generated by the $[g_i]$ and
$\ell[h_i].$

\begin{lemma}\label{lem4.4}
Keep the notation of the preceding paragraph, and suppose that the ramification index at all
primes over $\ell$ of $F$ is at most $\ell-1.$ Then $\ell$ exactly divides the order of
$(U_{1}(\mathfrak{n}_\Sigma).Z(\mathbb{A}_F^\infty)\cap h_i^{-1}D^\times h_i)/F^\times$.
\end{lemma}

\begin{Proof} One easily reduces the statement to showing that finite
subgroups of $D^\times$ having $\ell$-power order must have order exactly $1$ or $\ell$ (use the
two exact sequences in the proof of Lemma 1.1 of \cite{taylor}). Further, there can be a
non-trivial finite subgroup of $\ell$-power order if and only if $\zeta+\zeta^{-1}$ is in $F.$
Since any group of order $\ell^2$ is abelian, the only possible non-trivial finite subgroup has
to have order exactly $\ell$.
\end{Proof}

\begin{lemma}\label{lem4.5}
With the notation as above, let $f\in S_\mathcal{O}(U_1(\mathfrak{n}_\Sigma)).$ Then $T_x(f)\in
S_\mathcal{O}(U_1(\mathfrak{n}_\Sigma))^*$ for any prime $x\notin \Sigma$ with $\mathrm{\bf
N}x\equiv-1\pmod{\ell}$.
\end{lemma}

\begin{Proof}
Let $U^{(0)}$ be the subgroup of $U_{1}(\mathfrak{n}_\Sigma)$ consisting of elements whose $x$th
component is congruent to $
\begin{pmatrix} *&0\\ *& * \end{pmatrix}\pmod{\varpi_x}.$
Let $\zeta\in h^{-1}D^\times h\cap U_{1}(\mathfrak{n}_\Sigma).Z(\mathbb{A}_F^\infty)$ have order
exactly $\ell$ in the quotient $(h^{-1}D^\times h\cap
U_1(\mathfrak{n}_\Sigma).Z(\mathbb{A}_F^\infty))/F^\times$. We need to compute $T_x(f)(h)$ and
check that it is a multiple of $\ell$. Starting with a double coset decomposition given by
$\coprod_{i=0}^{\ell-1}\zeta^i\ *\ U^{(0)} $ and using the fact that $\zeta\notin U^{(0)},$ we
get a disjoint decomposition
$$U_{1}(\mathfrak{n}_\Sigma)=\coprod_{i=1}^{\ell}\coprod_{j=1}^{(\mathrm{\bf
N}x+1)/\ell}\zeta^iu_jU^{(0)}.$$ This shows that, by index considerations,
$$U_1(\mathfrak{n}_\Sigma)\begin{pmatrix}\varpi_x&0\\0&1\end{pmatrix}
U_1(\mathfrak{n}_\Sigma)=\coprod_{i=1}^{\ell}\coprod_{j=1}^{(\mathrm{\bf
N}x+1)/\ell}\zeta^iu_j\begin{pmatrix}\varpi_x&0\\0& 1\end{pmatrix}U_1(\mathfrak{n}_\Sigma).$$
Since $h\zeta^i=d_ih$ for some $d_i\in D^\times$, we have
\begin{eqnarray*}
T_x(f)(h) &=& \sum_{i=1}^{\ell}\sum_{j=1}^{(\mathrm{\bf
N}x+1)/\ell}f\left(h\zeta^iu_j\begin{pmatrix} \varpi_x & 0\\ 0&
1\end{pmatrix}\right)\\
&=& \sum_{i=1}^{\ell}\sum_{j=1}^{(\mathrm{\bf N}x+1)/\ell}f\left(hu_j\begin{pmatrix} \varpi_x &
0\\ 0&
1\end{pmatrix}\right)\\
&=&\ell\sum_{j=1}^{(\mathrm{\bf N}x+1)/\ell}f\left(hu_j\begin{pmatrix} \varpi_x & 0\\ 0&
1\end{pmatrix}\right).
\end{eqnarray*}
The lemma follows.
\end{Proof}

Now we discuss various properties of the modular forms and Hecke operators.
\begin{theorem} \label{th4.6}
Keeping the assumptions of the two preceding lemmas, we have the following:
\begin{enumerate}
\item The $\mathcal{O}$-module $S_\mathcal{O}(U_1(\mathfrak{n}_\Sigma))^*$ is invariant under the
action of Hecke operators. \item The pairing on $S_K(U_1(\mathfrak{n}_\Sigma))$ induces a perfect
pairing
$$S_\mathcal{O}(U_1(\mathfrak{n}_\Sigma))\times S_\mathcal{O}(U_1(\mathfrak{n}_\Sigma))^*
\longrightarrow\mathcal{O}.$$ \item Let $\mathfrak{m}$ be a non-Eisenstein maximal ideal of the
Hecke algebra $\mathbb{T}_\mathcal{O}(U_1(\mathfrak{n}_\Sigma)).$ Then
$S_\mathcal{O}(U_1(\mathfrak{n}_\Sigma))_\mathfrak{m}=
S_\mathcal{O}(U_1(\mathfrak{n}_\Sigma))^*_\mathfrak{m}.$ As a consequence, the pairing on
$S_K(U_1(\mathfrak{n}_\Sigma))$ induces a perfect pairing on
$S_\mathcal{O}(U_1(\mathfrak{n}_\Sigma))_\mathfrak{m}.$
\end{enumerate}
\end{theorem}

\begin{Proof}
The first part is easily checked using the given pairing on $S_K(U_1(\mathfrak{n}_\Sigma)).$ The
second part follows from Lemma~\ref{lem4.4}. The third part is a direct
consequence of Lemma~\ref{lem4.5}.
\end{Proof}

\section{Deformations in the minimal case}
In this section, we show that the universal deformation ring in the minimal case is isomorphic to
a Hecke algebra, and we show that these are complete intersection rings of relative dimension
zero over $\mathbb{Z}_p$.

Recall that we are given a continuous representation
$$\overline\rho:G_F\longrightarrow\mathrm{GL}_2(k)$$
satisfying the various properties listed in the beginning of section~\ref{sec3}, and also
satisfying Assumption~\ref{ass2}. In this and the next section, we shall assume the following
additional modularity condition.
\begin{assumption}\label{ass3}\rm
Let $U_0$ denote $U_{\{1\}}(\mathfrak{n}_\emptyset)$. Then we assume that there is a continuous
homomorphism $\phi :\mathbb{T}_\mathcal{O}(U_0)\rightarrow k$ with non-Eisenstein kernel which
gives our representation $\overline\rho.$ We write $\mathfrak{m}_\emptyset$ for the kernel.
\end{assumption}
Our aim is to show that the natural map
$R_\emptyset\twoheadrightarrow\mathbb{T}_\mathcal{O}(U_0)_{\mathfrak{m}_\emptyset}$ is an
isomorphism of complete intersection rings.

Fix a finite set of primes $\Sigma$ of $F$ not dividing $\ell$ such that for every $x\in\Sigma,$
we have
\begin{itemize}
\item $\mathrm{\bf N}x\equiv 1\pmod{\ell},$ \item $\overline{\rho}$ is unramified at $x$ and
 has distinct eigenvalues $\alpha_x\neq\beta_x.$
\end{itemize}
We denote the maximal $\ell$-power quotient of $(\mathcal{O}_F/x)^\times,$ for $x\in\Sigma,$ by
$\Delta_x$ and set $\Delta_\Sigma=\prod\Delta_x.$ We define the following objects (all products
are over $x\in\Sigma$):
\begin{enumerate}
\item an ideal $\mathfrak{n}_\Sigma=\prod x^2$. \item compact open subgroups
$U_{0,\Sigma}=U_{\{1\}}(\mathfrak{n}_\Sigma)$ and
$U_{1,\Sigma}=U_{\Delta_\Sigma}(\mathfrak{n}_\Sigma).$ \item an ideal $\mathfrak{m}_\Sigma$ of
either $\mathbb{T}(U_{0,\Sigma})$ or $\mathbb{T}(U_{1,\Sigma})$ generated by $\ell$ and
\begin{itemize}
\item $T_x-\mathrm{tr}\,\overline{\rho}(\mathrm{Frob}_x)$ for $x\nmid\ell\mathfrak{n}_\Sigma,$
and \item $\mathrm{\bf U}_{\varpi_x}-\alpha_x$ for $x\in\Sigma.$
\end{itemize}
\end{enumerate}

Note that Lemma 2.1 and Lemma 2.2 of \cite{taylor} remain true in the present situation (and we
will write them down again in a moment). We also have the fact that $S_\mathcal{O}(U_{1,\Sigma})$
is an $\mathcal{O}[\Delta_\Sigma]$-module via $h\rightarrow\langle h\rangle.$ But slight care is
required for the critical Lemma 2.3 and Corollary 2.4 of \cite{taylor}: it is no longer obvious
that $S_\mathcal{O}(U_{1,\Sigma})_{\mathfrak{m}_\Sigma}$ is free over
$\mathcal{O}[\Delta_\Sigma].$ Nonetheless, we can still get the `patching modules' technique of
\cite{Di} to work.

We first present a trivial reformulation of Theorem 2.1 of \cite{Di}.

\begin{theorem}\label{th5.2}
Fix a positive integer $r,$ a finite field $k$; set $A=k[[S_1,\ldots ,S_r]]$ and $B=k[[X_1,\ldots
, X_r]].$ We denote the maximal ideal of $A$ by $\mathfrak{n}.$ We are given: a $k$-algebra $R,$
a non-zero $R$-module $H$ which is finite dimensional over $k.$ For each positive integer $n,$ we
suppose that we have $k$-algebra homomorphisms $\phi_n:A\rightarrow B$ and $\psi_n:B\rightarrow
R,$ a $B$-module $H_n$ and a $B$-linear homomorphism $\pi_n:H_n\rightarrow H$ such that:
\begin{itemize}
\item $\psi_n$ is surjective and $\psi_n\phi_n=0$, \item $\pi_n$ induces an isomorphism between
$H_n/\mathfrak{n}H_n$ and $H,$ and \item there is an unbounded sequence of positive integers
$(a_n)_{n\geq 1}$ such that $H_n/\mathfrak{n}^{a_n}H_n$ is free over $A/\mathfrak{n}^{a_n}.$
\end{itemize}
Then $R$ is a complete intersection, and $H$ is free over $R.$
\end{theorem}

We now begin analyzing and comparing the $\mathcal{O}[\Delta_\Sigma]$-module structures of
$S_\mathcal{O}(U_{0,\Sigma})$ and $S_\mathcal{O}(U_{1,\Sigma}).$ Denote the augmentation ideal of
$\mathcal{O}[\Delta_\Sigma]$ by $I_{\Delta_\Sigma}.$ Obviously, functions in
$S_\mathcal{O}(U_{0,\Sigma})$ are precisely the elements of $S_\mathcal{O}(U_{1,\Sigma})$ which
are invariant under the action of $\Delta_\Sigma ;$ there is a `norm' map
$$\sum_{h\in\Delta_\Sigma}\langle h\rangle
:S_\mathcal{O}(U_{1,\Sigma})_{\Delta_\Sigma}\longrightarrow S_\mathcal{O}(U_{0,\Sigma}),$$ where
the subscript denotes coinvariants.

\begin{proposition} The norm map
 $$\sum_{h\in\Delta_\Sigma}\langle h\rangle
:S_\mathcal{O}(U_{1,\Sigma})\longrightarrow S_\mathcal{O}(U_{0,\Sigma})$$ has kernel
$I_{\Delta_\Sigma}S_\mathcal{O}(U_{1,\Sigma})$ and surjects onto $S_\mathcal{O}(U_{0,\Sigma})^*.$

The $\mathbb{T}(U_{1,\Sigma})$-module
$$\left(\sum_{h\in\Delta_\Sigma[\ell]}h\right)S_\mathcal{O}(U_{1,\Sigma})$$
is free over $\mathcal{O}[\Delta_\Sigma/\Delta_\Sigma[\ell]];$ and the norm map factorizes, in an
obvious way, as the composite of
$$\sum_{h\in\Delta_\Sigma[\ell]}\langle h \rangle \quad\mathrm{and}\quad\sum_{h\in\Delta_\Sigma/
\Delta_\Sigma[\ell]}\langle h\rangle.$$
\end{proposition}

\begin{Proof} We have a decomposition
$$(D\otimes_\mathbb{Q}\mathbb{A}^\infty)^\times=\coprod D^\times
t_iU_{0,\Sigma}.Z(\mathbb{A}_F^\infty).$$ For $h\in \Delta_\Sigma,$ we have a lift $\tilde{h}\in
(\mathbb{A}_F^\infty)^\times$ which gives the coset decomposition
$$U_{0,\Sigma}=\coprod_{h\in\Delta_\Sigma}\begin{pmatrix}\tilde{h}
&0\\ 0&1\end{pmatrix}U_{1,\Sigma}.$$ There is an obvious transitive action of $\Delta_\Sigma$ on
this coset decomposition.

For each $t_i$, we define
$$\mathrm{Stab}_i =\left\{ h\in\Delta_\Sigma\,\left|\,
D^\times t_iU_{1,\Sigma}.Z(\mathbb{A}_F^\infty)=D^\times
t_i\begin{pmatrix}\tilde{h} &0\\
0&1\end{pmatrix}U_{1,\Sigma}.Z(\mathbb{A}_F^\infty)\right.\right\}.$$ Obviously, the definition
is independent of the representatives $t_i$ and depends only the double coset decomposition. We
get the double coset decomposition
$$(D\otimes_\mathbb{Q}\mathbb{A}^\infty)^\times=\coprod_{i}\coprod_{h\in\Delta_\Sigma/\mathrm{Stab}_i} D^\times
t_i\begin{pmatrix}\tilde{h} &0\\
0&1\end{pmatrix}U_{1,\Sigma}.Z(\mathbb{A}_F^\infty).$$ In particular, we see that the set
$$\bigcup_{i}\left\{ \langle h\rangle
[t_i]\,|\,h\in\Delta_\Sigma/\mathrm{Stab}_i\right\}$$ is a basis for the free
$\mathcal{O}$-module $S_\mathcal{O}(U_{1,\Sigma}).$

It is now clear that the image of the map
$$\sum_{h\in\Delta_\Sigma}\langle h\rangle
:S_\mathcal{O}(U_{1,\Sigma})\longrightarrow S_\mathcal{O}(U_{0,\Sigma})$$ is free over
$\mathcal{O}$ with basis $\left\{|\mathrm{Stab}_i|[t_i]\right\}_i.$ The fact that the kernel is
the image of the augmentation ideal is obvious once we show that it is enough to consider
elements in the kernel having the form
$$x=\sum_{h\in\Delta_\Sigma/\mathrm{Stab}_i}a_h\langle h\rangle
[t_i]\quad \mathrm{with}\quad
a_h\in\mathcal{O}\quad\mathrm{and}\quad\sum_{h\in\Delta_\Sigma/\mathrm{Stab}_i}a_h=0.$$ It
suffices to consider such $x$ because we can write $x=\sum x_i$, where $x_i$ lies in the kernel
and has the form $|\mathrm{Stab}_i|(\sum a_h)[t_i].$

We now show that the image of the norm map is $S_\mathcal{O}(U_{0,\Sigma})^*$ by proving that the
order of $\textrm{Stab}_i$ is equal to the power of $\ell$ that divides the order of
$\left(t_i^{-1}D^\times t_i\cap U_{0,\Sigma}.Z(\mathbb{A}_F^\infty)\right)/F^\times .$

We claim that the order of  $\left(t_i^{-1}D^\times t_i\cap
U_{1,\Sigma}.Z(\mathbb{A}_F^\infty)\right)/F^\times $ is not divisible by $\ell.$ Indeed, let
$\alpha\in t_i^{-1}D^\times t_i\cap U_{1,\Sigma}.Z(\mathbb{A}_F^\infty)$ be such that
$\alpha^\ell \in F^\times.$ Fix a place $x\in \Sigma.$ We can write the $x$-th component of
$\alpha \in U_{1,\Sigma}.Z(\mathbb{A}_F^\infty)$ as $u_xz_x$ where $z_x\in K_x$ and $u_x\in
GL_2(\mathcal{O}_x)$ satisfies
$$u_x\equiv \begin{pmatrix} h& *\\
0&1\end{pmatrix}\pmod{\omega_x}$$ with $h$ having order prime to $\ell.$ Raising $u_x$ to the
$\ell-$th power, one deduces that $u_x$ reduces to the identity mod $\omega_x ,$ and hence that
$u_x$ is trivial. This then implies that $\alpha\in F^\times.$

Let $m$ be the prime to $\ell$ part of the order of $\left(t_i^{-1}D^\times t_i\cap
U_{0,\Sigma}.Z(\mathbb{A}_F^\infty)\right)/F^\times.$ We define a map
$\theta:\textrm{Stab}_i\longrightarrow  \left(t_i^{-1}D^\times t_i\cap
U_{0,\Sigma}.Z(\mathbb{A}_F^\infty)\right)/F^\times$ as follows: If $h\in \mathrm{Stab}_i,$ we
must have $t_i^{-1}dt_i=hu_1a=x$ (say) for some $d\in D^\times,\ u_1\in U_{1,\Sigma}$ and $a\in
(\mathbb{A}_F^\infty)^\times.$ Thus $x\in t_i^{-1}D^\times t_i\cap
U_{0,\Sigma}.Z(\mathbb{A}_F^\infty),$ and we set $\theta(h)=x^m\pmod{F^\times}.$ By the claim
established in the previous paragraph, it follows that $\theta$ is a well-defined injective
homomorphism from $\textrm{Stab}_i$ to the $\ell$-primary part of $\left(t_i^{-1}D^\times t_i\cap
U_{0,\Sigma}.Z(\mathbb{A}_F^\infty)\right)/F^\times.$ Since by Lemma~\ref{lem4.4} the order of
the $\ell$-primary part of $\left(t_i^{-1}D^\times t_i\cap
U_{0,\Sigma}.Z(\mathbb{A}_F^\infty)\right)/F^\times$ is exactly $\ell$ or 1, it is then simple to
verify that $\theta$ is an isomorphism between $\textrm{Stab}_i$ and the $\ell$-primary part of
$\left(t_i^{-1}D^\times t_i\cap U_{0,\Sigma}.Z(\mathbb{A}_F^\infty)\right)/F^\times.$ It follows
that the image of the norm map is exactly $S_\mathcal{O}(U_{0,\Sigma})^*.$

The last part of the proposition follows since $\mathrm{Stab}_i\subset\Delta_\Sigma[\ell].$
\end{Proof}

The following is Lemma 2.2 of \cite{taylor}. The proof given in \cite{taylor} works verbatim in
our case (thanks to Theorem~\ref{th4.6}).

\begin{lemma} There is an isomorphism
$S_\mathcal{O}(U_{0,\emptyset})_{\mathfrak{m}_\emptyset}\rightarrow
S_\mathcal{O}(U_{0,\Sigma})_{\mathfrak{m}_\Sigma}$ inducing an isomorphism
$\mathbb{T}(U_{0,\Sigma})_{\mathfrak{m}_\Sigma}\rightarrow
\mathbb{T}(U_{0,\emptyset})_{\mathfrak{m}_\emptyset} .$ $\hfill\Box$
\end{lemma}

Using the fact that the rings in consideration are semi-local, reduced and complete (they are
finite flat $\mathbb{Z}_\ell$-algebras), and Theorem~\ref{th4.6}, we get the following:
\begin{corollary}
\begin{enumerate}
\item There is an isomorphism
$S_\mathcal{O}(U_{1,\Sigma})_{\mathfrak{m}_\Sigma,\Delta_\Sigma}\longrightarrow
S_\mathcal{O}(U_{1,\emptyset})_{\mathfrak{m}_\emptyset}.$ This isomorphism is compatible with the
map on Hecke algebras $\mathbb{T}(U_{1,\Sigma})_{\mathfrak{m}_\Sigma}\rightarrow
\mathbb{T}(U_{0,\emptyset})_{\mathfrak{m}_\emptyset}$ which sends:
\begin{itemize}
\item $T_x$ to $T_x$ for $x$ not dividing $\ell\mathfrak{n}_\Sigma,$ \item $\langle h\rangle$ to
$1$ for $h\in \Delta_\Sigma,$ and \item $\mathrm{\bf U}_{\varpi_x}$ to $A_x$ for $x\in \Sigma $
where $A_x$ is the unique root of $X^2-T_xX+\mathrm{\bf N}x$ in
$\mathbb{T}(U_{0,\emptyset})_{\mathfrak{m}_\emptyset}$ congruent to
$\alpha_x\pmod{\mathfrak{m}_\emptyset}$.
\end{itemize}
\item The surjection $S_\mathcal{O}(U_{1,\Sigma})_{\mathfrak{m}_\Sigma }\twoheadrightarrow
S_\mathcal{O}(U_{1,\emptyset})_{\mathfrak{m}_\emptyset}$ given by composing the norm map with the
isomorphism of the preceding lemma factorizes as the composite of
$$S_\mathcal{O}(U_{1,\Sigma})_{\mathfrak{m}_\Sigma,\Delta_\Sigma}\twoheadrightarrow
H_\Sigma\quad\mathrm{and}\quad H_\Sigma\longrightarrow
S_\mathcal{O}(U_{1,\emptyset})_{\mathfrak{m}_\emptyset}$$ where:
\begin{itemize}\item $H_\Sigma$ is a
$\mathbb{T}(U_{1,\Sigma})_{\mathfrak{m}_\Sigma}$-algebra and the maps are compatible with the
algebra structures, and \item $H_\Sigma$ is a free
$\mathcal{O}[\Delta_\Sigma/\Delta_\Sigma[\ell]]$ module.
\end{itemize}
\end{enumerate}
$\hfill\Box$
\end{corollary}

We apply the above corollary to the sets $\Sigma_n$ produced by Theorem~\ref{th3.2-2}. Applying
the `patching modules' result of Diamond~\cite{Di} and Fujiwara~\cite{Fu} (Theorem~\ref{th5.2}
above), we get the following result.

\begin{theorem}
 The natural map
$$R_\emptyset\longrightarrow\mathbb{T}(U_0)_{\mathfrak{m}_\emptyset}$$
is an isomorphism of complete intersection rings and the module
$S_\mathcal{O}(U_0)_{\mathfrak{m}_\emptyset}$ is free over
$\mathbb{T}(U_0)_{\mathfrak{m}_\emptyset}.$ $\hfill\Box$
\end{theorem}

\section{Non-minimal level}

The proof of the result in the non-minimal case given in \cite{taylor} remains valid in our case.
We shall only give a sketch. Throughout this section, we keep the various assumptions (and
notation) of the last section.

Fix a homomorphism $\pi_\emptyset : R_\emptyset\twoheadrightarrow\mathcal{O}.$ We now let
$\Sigma$ be a finite set of primes of $F$ not containing any primes above $\ell.$ We denote by
$\pi_\Sigma$ the surjection $R_\Sigma\twoheadrightarrow\mathcal{O}$ obtained by taking the
composite of
$$R_\Sigma\twoheadrightarrow R_\emptyset\twoheadrightarrow\mathcal{O}$$
where the first map is the one given by the universal property of $R_\Sigma$ and the second map
is $\pi_\emptyset.$ We shall denote the kernel of $\pi_\Sigma$ by $\mathfrak{P}_\Sigma .$

Let $\mathfrak{n}_\Sigma=\prod_{x\in\Sigma}x^2,$ and let
$U_\Sigma=U_{\{1\}}(\mathfrak{n}_\Sigma).$ Also, let $\mathfrak{m}_\Sigma$ be the maximal ideal
of $\mathbb{T}_\mathcal{O}(U_\Sigma)$ corresponding to our residual representation $\overline\rho
.$ We denote by $\mathbb{T}_\Sigma$ the localization
$\mathbb{T}_\mathcal{O}(U_\Sigma)_{\mathfrak{m}_\Sigma},$ and write $S_\Sigma$ for the
$\mathbb{T}_\Sigma$-module $S_\mathcal{O}(U_\Sigma)_{\mathfrak{m}_\Sigma}.$

We then have the following.

\begin{theorem} The natural map
$R_\Sigma\twoheadrightarrow\mathbb{T}_\Sigma$ is an isomorphism of complete intersection rings
and $S_\Sigma$ is free over $\mathbb{T}_\Sigma.$
\end{theorem}

To prove the theorem, one needs to check (by Theorem 2.4 of \cite{Di}) that the order of
$\mathfrak{P}_\Sigma/\mathfrak{P}_\Sigma^2$ divides the order of
$$\Omega_\Sigma\stackrel{\mathrm{def}}{=}\frac{S_\Sigma}{S_\Sigma[\mathfrak{P}]\oplus
S_\Sigma[\mathrm{Ann}_{\mathbb{T}_\Sigma}\mathfrak{P}]}.$$ A standard computation shows that the
order of $\mathfrak{P}_\Sigma /\mathfrak{P}_\Sigma^2$ divides
$$\#\left(\mathfrak{P}_\emptyset/\mathfrak{P}_\emptyset^2\right)
\prod_{x\in\Sigma}\#\left(\mathcal{O}/ (1-\mathrm{\bf N}x)(T_x^2-(1+\mathrm{\bf
N}x)^2)\mathcal{O}\right),$$ and we shall prove that this expression is the order of
$\Omega_\Sigma$.

Note that $S_\Sigma[\mathfrak{P}_\Sigma]$ is a free $\mathcal{O}$-module of rank 1. Fix a perfect
symmetric $\mathcal{O}$-valued $\mathcal{O}$-bilinear pairing $\{\ ,\ \}_\Sigma$ on
$S_\Sigma[\mathfrak{P}_\Sigma],$ and let $j_\Sigma:S_\Sigma[\mathfrak{P}_\Sigma]\hookrightarrow
S_\Sigma$ be the natural inclusion. Also, define a pairing $\langle \ ,\ \rangle_\Sigma$ on
$S_\Sigma$ by
$$\langle f_1, f_2\rangle_\Sigma =(f_1, w_\Sigma f_2)$$ where $(\
,\ )$ is the standard pairing, and $w_\Sigma\in
\mathrm{GL}_2(\mathbb{A}_F^\infty)\cong(D\otimes_\mathbb{Q}\mathbb{A}^\infty)^\times$ is the
element defined by
$$ w_{\Sigma ,x} =\left\{ \begin{matrix}\mathrm{identity}, &\mathrm{if}\ x\notin\Sigma ,\\
 \begin{pmatrix} 0& 1\\ \varpi_x^2 & 0\end{pmatrix},& \mathrm{if}\
x\in\Sigma.\end{matrix}\right.$$ This new pairing is perfect, and the Hecke operators are
self-adjoint with respect to $\langle \ ,\ \rangle_\Sigma .$

Now let $x$ be a prime not dividing $\mathfrak{n}_\Sigma\ell.$ There is a well-defined map
$$i_x:S_\Sigma\longrightarrow
S_{\Sigma\cup\{x\}}$$ which is obtained from the map sending $f\in S_\mathcal{O}(U_\Sigma)$ to
$$(\mathrm{\bf N}x)f-\begin{pmatrix}1&0\\
0&\varpi_x\end{pmatrix} T_xf+\begin{pmatrix}1&0\\
0&\varpi_x^2\end{pmatrix}f\ \in S_\mathcal{O}(U_{\Sigma\cup\{x\}}).$$ Under this map, the image
of $S_\Sigma[\mathfrak{P}_\Sigma]$ is contained in
$S_{\Sigma\cup\{x\}}[\mathfrak{P}_{\Sigma\cup\{x\}}].$ We denote by $\widetilde{i_x}$ the
resulting map from $S_\Sigma[\mathfrak{P}_\Sigma]$ to
$S_{\Sigma\cup\{x\}}[\mathfrak{P}_{\Sigma\cup\{x\}}].$

We then have the following.
\begin{itemize}
\item Let $i_x^*$ be the adjoint of $i_x$ with respect to the pairings $\langle \ ,\
\rangle_\Sigma$ and $\langle \ ,\ \rangle_{\Sigma\cup\{x\}}.$ The composite $i_x^*\circ i_x$ is
equal to $$\mathrm{\bf N}x(1-\mathrm{\bf N}x)(T_x^2-(1+\mathrm{\bf N}x)^2).$$ \item
$i_x(S_\Sigma[\mathfrak{P}_\Sigma])=S_{\Sigma\cup\{x\}}[\mathfrak{P}_{\Sigma\cup\{x\}}]$. This
follows from Ihara's lemma (see Lemma 3.1 of \cite{taylor}). \item Let $j_\Sigma^*$ be the
adjoint of $j_\Sigma$ with respect to the pairings $\{\ ,\ \}_\Sigma$ and $\langle \ ,\
\rangle_\Sigma .$ It induces an isomorphism
$$j_\Sigma^*:\Omega_\Sigma\stackrel{\sim}{\longrightarrow}
\frac{S_\Sigma[\mathfrak{P}_\Sigma]}{j_\Sigma^*S_\Sigma[\mathfrak{P}_\Sigma]}.$$ \item Let
$\widetilde{i_x}^*$ be the adjoint of $\widetilde{i_x}$ with respect to the pairings $\{ \ ,\
\}_\Sigma$ and $\{ \ ,\ \}_{\Sigma\cup\{x\}}.$ It is an isomorphism, and we have
$\widetilde{i_x}^*\circ j_{\Sigma\cup\{x\}}^*=j_\Sigma^*\circ i_x^*.$
\end{itemize}

It follows that
$$\#\Omega_\Sigma=\#\Omega_\emptyset\prod_{x\in\Sigma}\#\left(\mathcal{O}/
(1-\mathrm{\bf N}x)(T_x^2-(1+\mathrm{\bf N}x)^2)\mathcal{O}\right).$$ The result in the minimal
case implies that $\#\Omega_\emptyset=\#(\mathfrak{P}_\emptyset/\mathfrak{P}_\emptyset^2),$ and
hence that
$$\left.\#\frac{\mathfrak{P}_\Sigma}{\mathfrak{P}_\Sigma^2}\right|\#\Omega_\Sigma.$$

\section{Modularity of Galois representations and elliptic curves}
\label{modularity}

We now collect the results of the preceding two sections.

Let $F$ be a totally real, finite extension of $\mathbb{Q}.$ Let $\mathcal{O}$ be the ring of
integers in a finite extension of $\mathbb{Q}_\ell$ where $\ell$ is an odd prime, and let $k$ be
its residue field. We suppose that we are given continuous representations
$$\rho_i:G_F\longrightarrow\mathrm{GL}_2(\mathcal{O}),\quad i=1,2$$
satisfying the following properties:
\begin{itemize}
\item $\rho_i$ ($i=1,2$) is an odd representation unramified outside finitely many primes; \item
$\mathrm{det}\,\rho_1=\mathrm{det}\,\rho_2=\epsilon_\ell$ where $\epsilon_\ell$ is the
$\ell$-adic cyclotomic character. \item The residual representations
$\overline\rho_i:G_F\rightarrow \mathrm{GL}_2(k)$ are equivalent and are absolutely irreducible.
We denote the residual representation by $\overline\rho.$
\end{itemize}

\begin{theorem}\label{th7.1}
With notations as in the preceding paragraph, we make the following assumptions.
\begin{itemize}
\item The restriction of $\overline\rho$ to the absolute Galois group of $F(\zeta_\ell)$ is
absolutely irreducible; furthermore,  if $\ell=5$ and $\textrm{Proj}\
\overline\rho|_{\mathrm{Gal}\left(\overline{F}/F\left(\zeta_\ell \right)\right)}$ then
$[F(\zeta_\ell):F]=4.$ \item (Conditions at $\ell.$) Let $v$ be any prime of $F$ dividing $\ell
,$ and let $I_v$ be the inertia group of $F_v.$ We assume:
\begin{enumerate}
\item $\overline\rho|_{I_v}\sim\Omega_2|_{I_v}$, where $\Omega_2$ is the second fundamental
character of the inertia group of $\mathbb{Q}_\ell.$
\item Let $\mathfrak{m}$ be the maximal ideal of $\mathcal{O},$ and let $\overline\rho_{i,n}$ be
the reduction of $\rho_i$ modulo $\mathfrak{m}^n.$ Then $\overline\rho_{i,n}|_{F_v}$ is finite
flat.
\end{enumerate}
\item The ramification index of $F$ at any prime above $\ell$ is less than or equal to $\ell-1.$
\end{itemize}
Under these assumptions, the modularity of $\rho_1$ implies the modularity of $\rho_2.$
\end{theorem}

\begin{Proof}
We can find a totally real, finite soluble extension $F'/F$ such that:
\begin{itemize}
\item The extension $F'/F$ is unramified at primes dividing $\ell.$ \item
$\overline\rho|_{G_{F'}}$ satisfies Assumption~\ref{ass3}. (For this, we need to use the
modularity of $\rho_1$ along with the base change results in \cite{SW3}.)
\end{itemize}
It follows that $\rho_2|_{G_{F'}}$ is modular. Langlands' cyclic base change then shows that
$\rho_2$ is modular. \end{Proof}

In section~\ref{AppII}, we will give some applications to the modularity of elliptic curves.
However, let us remark here that Theorem~\ref{th7.1} will not apply in general to all
supersingular curves, as the first condition at $\ell$ will not be satisfied in general. Indeed,
let $F=\mathbb{Q}(\sqrt{3})$, and let $E$ denote the elliptic curve
$$y^2=x^3+\sqrt{3}x^2+x+1.$$
The curve has discriminant $32(3\sqrt{3}-14)$, and hence has good reduction at the prime
$\sqrt{3}$ above 3. On the other hand, it is easy to show that multiplication by 3 on the group
law of an elliptic curve
$$y^2=x^3+a_2x^2+a_4x+a_6$$
is given by
$$[3]t=3t-8a_2t^3+\cdots,$$
so that the curve above has supersingular reduction at $\sqrt{3}$, as $v_3(a_2)=v_3(\sqrt{3})>0$,
showing that the formal group at 3 has height 2. As in Serre~\cite{serre}, Proposition 10, the
action of tame inertia on the 3-torsion points is given by 2 copies of the fundamental character
of level 1, rather than by the fundamental character of level 2.

Serre's argument also shows that in order that the mod 3 representation of the curve $E$ be given
(on tame inertia) by the fundamental character of level 2, it is necessary and sufficient that
the Newton polygon of the multiplication-by-3 map on the formal group should consist of a single
line from $(1,e)$ to $(9,0)$. This is automatic when $e=1$, but if $e>1$, then other situations
may arise, as above.

It follows that our main result can apply to all supersingular curves defined over fields $F$
unramified at 3, as well as to many examples of curves defined over more general fields.

\section{Applications I}

\begin{theorem}
Let $\overline\rho:\mathrm{Gal}(\overline{\mathbb{Q}}/\mathbb{Q})\longrightarrow
GL_2(\mathbb{F}_7)$ be an absolutely irreducible, continuous, odd representation. If the
projective image of $\overline\rho$ is insoluble, we also assume that:
\begin{itemize}
\item The projective image of inertia at $3$ has odd order. \item The determinant of
$\overline\rho$ restricted to the inertia group at $7$ has even order.
\end{itemize}
Then $\overline\rho$ is modular.
\end{theorem}

\noindent {\em Sketch of proof.} Of course, we need  only consider the case when the image of
$\overline\rho$ is insoluble. Moreover  by \cite{mano}, we can assume that the restriction of
$\overline\rho$ to a decomposition group at 7 is irreducible. Twisting by a quadratic character,
we can also assume that $\overline\rho|_{I_7}$ is equivalent to $\omega_2\oplus\omega_2^{7}$ or
$\omega_2^{13}\oplus\omega_2^{7.13}$ where $\omega_2:I_7\longrightarrow \mathbb{F}_{49}^\times$
is the second fundamental character. Applying the axiomatic formulation of Ramakrishna's result
in \cite{taylorartin}, together with  Theorems 3.2.1, 4.2.1 of \cite{conrad2}, one deduces the
existence of a  continuous, odd representation
$$\rho:\mathrm{Gal}(\overline{\mathbb{Q}}/\mathbb{Q})\longrightarrow
GL_2(\mathbb{Z}_7)$$ lifting $\overline\rho,$ unramified outside finitely many primes,
determinant the cyclotomic character times a finite order character, and  such that the Artinian
quotients $\rho\pmod{7^n}$ are finite flat when restricted to the absolute Galois group of
$\mathbb{Q}_7(7^{1/4}).$ Assuming the existence of a totally real soluble extension
$F/\mathbb{Q}$ such that $\overline\rho|_{G_F}$ is modular and the ramification index of
$F/\mathbb{Q}$ at 7 is at most 6, one deduces the modularity of $\rho$ by Theorem~\ref{th7.1} and
Langlands' cyclic base change.

We now explain how to find such a field $F.$ Firstly, we can find a finite soluble, totally real
extension $F_1/\mathbb{Q}$ and a quadratic twist of $\overline\rho|_{G_{F_1}},$ which we denote
by $\widetilde{\rho},$ such that the following conditions are satisfied.
\begin{itemize}
\item The determinant of $\widetilde{\rho}$  is  the mod $7$ cyclotomic character. \item
Conditions at $3$: Let $v$ be any prime of $F_1$ above $3,$ and let $D_v$ be a decomposition
group at $v.$
\begin{itemize}
\item $\widetilde{\rho}$ is trivial on $D_v$. \item The ramification index of
$F_{1,v}/\mathbb{Q}_3$ is odd.
\end{itemize}
\item Conditions at $7$: Let $v$ be any prime of $F_1$ above $7,$ and let $D_v$, $I_v$ be the
decomposition and inertia groups at $v.$ Then, the ramification index of $F_{1,v}/\mathbb{Q}_7$
is exactly $4$. Furthermore, we have $\widetilde{\rho}|_{I_{F_{1,v}}}\cong
(\omega_2\oplus\omega_2^{7})|_{I_{F_{1,v}}}.$
\end{itemize}
 We denote by $X(\widetilde{\rho})$ the (completed) moduli space of elliptic curves with mod $7$
representation symplectically isomorphic to $\widetilde{\rho}$ (see \cite{mano} for details). The
canonical divisor embeds $X(\widetilde{\rho})$ as a quartic curve in $\mathbb{P}^2_{/F_1}.$

For each prime $v$ of $F_1$ dividing $3\infty,$ we can find a finite unramified  extension
$F_v/F_{1,v}$ and a line $L_v$ defined over $F_{1,v}$ such that  $L_v$ cuts
$X(\widetilde{\rho})_{/F_v}$ at four distinct points all of which are defined over $F_v.$
Moreover, the elliptic curves corresponding to these four points all have good ordinary reduction
when $v|3.$ (See  the fourth paragraph in section $5$ of \cite{mano}.) For primes above 7, we
have the following lemma:

\begin{lemma}
Let $v$ be a prime of $F_1$ above $7$. We can find a finite Galois extension $F_v/F_{1,v}$ and an
$F_v$-rational line $L_v$ such that the following holds.
\begin{itemize}
\item $L_v$ cuts $X(\widetilde{\rho})_{/F_v}$ at four distinct points all of which are defined
over $F_v.$ \item The ramification index of $F_v/\mathbb{Q}_7$ is at most $4.$ The four points of
intersection are all elliptic curves with good supersingular reduction.
\end{itemize}
\end{lemma}

Assuming the above lemma,  intersecting $X(\widetilde{\rho})$ with a line over $F_1$ which is
$v$-adically close to $L_v$ for each $v|3.7.\infty$ gives the following: There is a finite,
soluble, totally real $F\supset F_1\supset\mathbb{Q},$ and an elliptic curve $E_{/F}$ satisfying
the following conditions.
\begin{itemize}
\item $\overline\rho_{E,7}\sim\widetilde{\rho}|_{G_F}$ and
$\overline\rho_{E,3}:G_F\twoheadrightarrow GL_2(\mathbb{F}_3)$ is surjective. \item Conditions at
primes $v$ dividing $3$: $E$ has good ordinary reduction at every prime above $3$ and the
ramification index of $F$ at $3$ is odd. \item Conditions above $7$: $F/F_1$ is unramified at
every prime above $7$ and $E$ has good supersingular reduction at every prime above $7$.
\end{itemize}
The elliptic curve $E$ is modular by a result of Skinner and Wiles (\cite{SW2}), and therefore
$\overline\rho$ is also modular.$\hfill\Box$

\bigskip
\noindent {\em Proof of Lemma 8.2.} The modular curve
$X(\omega_2\oplus\omega_2^7)_{/\mathbb{Q}_7^\mathrm{nr}}$ is isomorphic to $X(\widetilde{\rho})$
over $\mathbb{Q}_7^\mathrm{nr}(\sqrt[4]{7}).$ The elliptic curve $ y^2=x^3+x$ has $j$-invariant
$1728$ and so has supersingular reduction. Taking a cyclic degree $3$ isogeny of $E$ if
necessary, we can assume that $X(\omega_2\oplus\omega_2^7)(\mathbb{Q}_7^{\mathrm{nr}})$ contains
an elliptic curve $E$ having good supersingular reduction and with $j$-invariant $1728.$ Let us
denote this point by $P.$ From the geometry of the Klein quartic (see the proposition in section
2 of \cite{elkies}), we see that there is a unique involution (in the automorphism group) fixing
$P.$ The normalizer of this involution is a Sylow 2-subgroup, and the orbit of $P$ when acted on
by the normalizer has size exactly 4. Furthermore, they (the points in the orbit) lie on a unique
line.

We can thus find a unique line $L$ passing through $P$ such that:
\begin{itemize}
\item $L$ is defined over $\mathbb{Q}_7^{\mathrm{nr}},$ \item $L$ passes through four distinct
points of $X(\omega_2\oplus\omega_2^7)$
 whose $j$-invariants
are $1728$.
\end{itemize}
We claim that two of these points are already defined over $\mathbb{Q}_7^{\mathrm{nr}}.$ We have
the point $P$ with corresponding elliptic curve $E.$ Note that $E$ has complex multiplication by
$\mathbb{Z}[i]$ (and the endomorphism ring is already defined over $\mathbb{Q}_7^{\mathrm{nr}}$).
We now check that the isogeny $E\stackrel{2-2i}{\longrightarrow}E$ gives us another point of
intersection (which is obviously defined over $\mathbb{Q}_7^{\mathrm{nr}}$). This can be checked
over $\mathbb{C},$ and follows from the following observations.
\begin{itemize}
\item The involution $\begin{pmatrix} 0 & 1\\-1&0\end{pmatrix}\in PSL_2(\mathbb{F}_7)$ fixes
$$(\{1/7, i/7\},\mathbb{C}/\mathbb{Z}+i\mathbb{Z})\in
X(7)(\mathbb{C}).$$ \item $\begin{pmatrix} 2 & 2\\-2&2\end{pmatrix}$ is in the normalizer of
$\begin{pmatrix} 0 & 1\\-1&0\end{pmatrix}$ and sends
$$(\{1/7,
i/7\},\mathbb{C}/\mathbb{Z}+i\mathbb{Z})\quad\mathrm{to}\quad (\{2-2i/7,
2+2i/7\},\mathbb{C}/\mathbb{Z}+i\mathbb{Z}).$$
\end{itemize}

Thus each of the four points of intersection are defined over
$\mathbb{Q}_7^\mathrm{nr}(\sqrt{7}).$ The Sylow 2-subgroup which acts transitively on these four
points is dihedral; in terms of generators and relations, it is given by
$$\langle\left. \alpha ,\beta \right|\alpha^4=\beta^2=e,
\beta\alpha\beta=\alpha^3\rangle.$$ The unique involution which stabilizes $P$ is $\alpha^2$, and
it is defined over $\mathbb{Q}_7^\mathrm{nr}.$ The other three points are given by $\alpha(P),
\beta(P)$ and $\alpha\beta(P)$.

We now check that $\alpha, \beta$ are defined over $\mathbb{Q}_7^\mathrm{nr}(\sqrt[4]{7}).$ If
$\sigma\in G_{\mathbb{Q}_7^\mathrm{nr}(\sqrt{7})},$ we have
$$(\sigma*\beta)(P)=\sigma(\beta(\sigma^{-1}P))=\beta(P).$$
Therefore, we have $\sigma*\beta=\alpha^{2i(\sigma)}\beta$ where
$$i:G_{\mathbb{Q}_7^\mathrm{nr}(\sqrt{7})}\longrightarrow
\mathbb{Z}/2\mathbb{Z}$$ is a continuous homomorphism which necessarily factors through
$\mathbb{Q}_7^\mathrm{nr}(\sqrt[4]{7}).$ Similarly for $\alpha .$ We can thus conclude that all
the four points of intersection have good supersingular reduction
$\mathbb{Q}_7^\mathrm{nr}(\sqrt[4]{7}).$

Finally, it follows that we can find a line defined over an extension of $F_{1,v}$ with absolute
ramification index 4 which cuts $X(\widetilde{\rho})$ at four distinct supersingular points, all
defined over that extension. Take $F_v$ to be the Galois closure of the extension thus
constructed, and take $L_v$ to be the line $L_{/F_v}.$ $\hfill\Box$

\section{Applications II}\label{AppII}
The aim of this section is to study the modularity of elliptic curves over certain totally real
fields, using Theorem~\ref{th7.1}. Our main results are given by Propositions~\ref{prop9.3}
and~\ref{prop9.4}. For the particular example of the field $\mathbb{Q}(\sqrt{2})$, we can prove
more; the analogue of the switch between $p=3$ and $p=5$ used by Wiles (\cite{Wi},~\S5) holds,
and we can use existing results, together with the new results in this paper, to deduce the
modularity of all semistable elliptic curves over $\mathbb{Q}(\sqrt{2})$.

In~\cite{JM}, it is explained that this implies a version of Fermat's Last Theorem over
$\mathbb{Q}(\sqrt{2})$. Further calculations in~\cite{JM} show that $\mathbb{Q}(\sqrt{2})$ is the
only real quadratic field for which one can hope to generalise the methods of Ribet and Wiles to
prove such a result. It seems remarkable to us that there are any fields other than $\mathbb{Q}$
for which all the numerology allows us to prove generalisations of Fermat's Last Theorem.

We begin by proving results for more general fields. We start with a preliminary lemma.

\begin{lemma}\label{lem9.2}
Let $p$ be equal to $3$ or $5$, and let $F$ be a totally real number field in which $p$ is
unramified. Let $E$ be an elliptic curve over $F$ with good supersingular reduction at some place
$v|p.$ Then
$$\left.\overline\rho_{E,p}\right|_{\mathrm{Gal}
\left(\overline{F}/F\left(\sqrt{(-1)^{(p-1)/2}p}\right)\right)}$$ is absolutely irreducible.
\end{lemma}

\begin{Proof} The presence of a non-trivial complex conjugation
shows that irreducibility is the same as absolute irreducibility for odd
$\mathrm{GL}_2(\mathbb{F}_p)$-valued representations of totally real fields. The lemma then
follows easily when $p=5.$

We now do $p=3.$ Suppose, for a contradiction, that the conclusion of the lemma fails. Let $I_v$
be a decomposition group at $v.$ Since the image $\overline\rho_{E,3}(I_v)$ is cyclic of order
$8,$ it follows that the image $\overline\rho_{E,3}(\mathrm{Gal}(\overline{F}/F))$ is the full
Sylow 2-subgroup of $\mathrm{GL}_2(\mathbb{F}_3).$ Denoting by $K$ the splitting field of
$\overline\rho_{E,3},$ it follows that the image
$\overline\rho_{E,3}(\mathrm{Gal}(K/F(\sqrt{-3})))$ is an abelian group of order $8.$

The Sylow 2-subgroup of $\mathrm{GL}_2(\mathbb{F}_3)$ is the group
$$\langle \left. c,\tau\right| c^2=\tau^8=1,\ c\tau=\tau^3
c\rangle ,$$ and we may suppose that
$$c=\begin{pmatrix}1&0\\ 0&-1\end{pmatrix},\quad
\tau=\begin{pmatrix}1&1\\ -1&1\end{pmatrix}.$$ Since the image of $\mathrm{Gal}(K/F(\sqrt{-3}))$
is in $\mathrm{SL}_2(\mathbb{F}_3)$, it must in fact be the subgroup generated by $\tau^2$ and
$c\tau.$ This subgroup is non-abelian, giving the desired contradiction.
\end{Proof}

The next two propositions prove modularity of many elliptic curves over certain totally real
fields, using Theorem~\ref{th7.1}.
\begin{proposition}\label{prop9.3}
Let $F$ be a totally real number field in which $3$ is unramified, and let $E$ be an elliptic
curve over $F$ with good supersingular reduction at primes above $3.$ Then $E$ is modular.
\end{proposition}

\begin{Proof} We proceed in several steps. By the result of
Langlands and Tunnell, we know that $\overline\rho_{E,3}$ is modular. However, in order to apply
Theorem~\ref{th7.1} we need to produce a modular lift with level coprime to $3$.

\noindent {\em Step I:} By Langlands' cyclic base change, we need only prove the result over a
totally real soluble extension. In particular, making an appropriate base change if necessary, we
can assume that $\overline\rho_{E,3}|_{D_v}$ is trivial for any prime $v|5.$

\noindent {\em Step II:} We can find an elliptic curve $E'$ over $F$ such that
\begin{itemize}
\item $\overline\rho_{E,3}\sim\overline\rho_{E',3}$, \item $\overline\rho_{E',5}$ has insoluble
image, \item $E'$ has good ordinary reduction at every prime above $5$
and $$\overline\rho_{E',5}|_{D_v}\cong \begin{pmatrix} * & *\\
0& *\end{pmatrix}\quad\mathrm{for~any}\quad v|5$$ with distinct characters on the diagonal, \item
$E'$ has good reduction at primes above $3$.
\end{itemize}
If we can show that $E'$ is modular, then $\rho_{E',3}$ will be a modular lift of
$\overline\rho_{E,3}$ of the `right level'; we can then use Theorem~\ref{th7.1} to conclude that
$\rho_{E,3}$ is modular.

In order to show that $E'$ is modular, we want to make use of its $5$-adic representation and
apply the results in \cite{SW2}. For this, we need to produce a nearly ordinary modular lift of
$\overline\rho_{E',5}.$ Again, we can work over totally real soluble extensions.

\noindent {\em Step III:} We can assume that $\overline\rho_{E',5}$ is trivial when we restrict
to primes above $3$. We can then find a second elliptic curve $E''$ such that
\begin{itemize}
\item $\overline\rho_{E',5}\sim\overline\rho_{E'',5}$, \item
$\overline\rho_{E'',3}:G_F\longrightarrow\mathrm{GL}_2(\mathbb{F}_3)$ is surjective, \item $E''$
has split multiplicative reduction at every prime above $3$
and $$\overline\rho_{E'',3}|_{D_v}\cong \begin{pmatrix} * & *\\
0& *\end{pmatrix}\quad\mathrm{for~any}\quad v|3$$ with distinct characters on the diagonal, \item
$E''$ has good ordinary reduction at primes above $5$.
\end{itemize}
By Theorem~\ref{th7.1}, $E''$ is modular.

Since $\rho_{E'',5}$ is a nearly ordinary modular lift, it follows that $\rho_{E',5}$ is modular.
\end{Proof}

\begin{proposition}\label{prop9.4}
Let $F$ be a totally real number field in which $3$ and $5$ are unramified. Let $E$ be an
elliptic curve over $F$ with semistable reduction at primes above $3$ and $5$. Further, assume
that $E$ has good supersingular reduction at primes above $5$ and that
$\overline\rho_{E,5}|_{\mathrm{Gal}(\overline{F}/F(\sqrt{5}))}$ is irreducible. Then $E$ is
modular.
\end{proposition}

\begin{Proof} Going up to a soluble totally real field
(without changing ramification at $3$ and $5$) if necessary, we can assume that
$\overline\rho_{E,5}|_{D_v}$ is trivial for places $v|3$ where $E$ has good reduction. Then using
the twisted modular curve $X(E[5])_{/F},$ we can find an elliptic curve $E'/F$ such that
\begin{itemize}
\item $\overline\rho_{E',5}\sim\overline\rho_{E,5}$, \item $E'$ has the same reduction type as
$E$ at primes above $5$, \item $E'$ is a Tate curve at primes above $3$, and \item
$\overline\rho_{E',3}: G_F\longrightarrow\mathrm{GL}_2(\mathbb{F}_3)$ is surjective.
\end{itemize}
It follows that $\rho_{E',3}$ is modular, and $\rho_{E',5}$ is a modular lift of
$\overline\rho_{E,5}$ of the `right level'. Therefore, using either Theorem~5.1 of \cite{SW2} or
Theorem~\ref{th7.1} of this article, it follows that $\rho_{E,5}$ is modular.
\end{Proof}

Having proven some results over general fields, we now specialise to the case
$F=\mathbb{Q}(\sqrt{2})$, for which, as we shall see, there is also a version of the switch
between 3 and 5 used by Wiles (\cite{Wi}, \S5). In particular, this allows us to prove the
modularity of all semistable elliptic curves over $\mathbb{Q}(\sqrt{2})$.

\begin{proposition}\label{prop9.5}
Let $E$ be a semistable elliptic curve over $\mathbb{Q}(\sqrt{2}).$ Let $p$ be either $3$ or $5$.
If $\overline\rho_{E,p}$ is irreducible, then
$$\left.\overline\rho_{E,p}\right|_{\mathrm{Gal}
\left(\overline{F}/F\left(\sqrt{(-1)^{(p-1)/2}p}\right)\right)}$$ is absolutely irreducible.
\end{proposition}

\begin{Proof} Suppose the proposition fails to hold. Then $p$ does
not divide the order of $\overline\rho_{E,p}(\mathrm{Gal}(\overline{F}/F)),$ and so the
semistability condition implies that $\overline\rho_{E,p}$ is unramified at primes not dividing
$p.$ Further, by Lemma~\ref{lem9.2}, we see that $E$ has good ordinary or multiplicative
reduction at $p$. Therefore, we must have
$$\left.\overline\rho_{E,p}\right|_{I_p}\sim\begin{pmatrix}\overline{\epsilon}_p
&0\\ 0&1\end{pmatrix}$$ where $\overline{\epsilon}_p$ is the mod $p$ cyclotomic character. (Note
also that $3$ and $5$ are inert in $\mathbb{Q}(\sqrt{2})$.)

Let $K$ be the splitting field of $\overline\rho_{E,p},$ and let $\zeta_p$ be a primitive $p$th
root of unity. Then $K$ is an everywhere unramified abelian extension of
$\mathbb{Q}(\sqrt{2},\zeta_p).$ The class number of $\mathbb{Q}(\sqrt{2},\zeta_p)$ is then
checked to be equal to $1$ for both $p=3$ and $p=5$ (we used {\tt PARI} to verify this), giving
the required contradiction.
\end{Proof}

\begin{proposition}\label{prop9.6}
The modular curve $X_0(15)$ has exactly eight $\mathbb{Q}(\sqrt{2})$-rational points. Four of
these are cusps. The remaining four are elliptic curves with additive reduction at $5$.
\end{proposition}

\begin{Proof}
$X_0(15)$ is an elliptic curve, and, using Cremona's tables~\cite{cremona}, we can find an
explicit equation for it. The rank of $X_0(15)$ regarded as an elliptic curve over
$\mathbb{Q}(\sqrt{2})$ is the sum of its rank over $\mathbb{Q}$ and the rank (over $\mathbb{Q}$)
of its quadratic twist. An equation of $X_0(15)$ over $\mathbb{Q}$ is $y^2+xy+y=x^3+x^2-10x-10$,
and its quadratic twist over $\mathbb(\sqrt{2})$ is $y^2=x^3+x^2-641x- 3105$, which is curve
960G3 in Cremona's tables. Both curves have rank 0 over $\mathbb{Q}$, and it follows that
$X_0(15)$ has rank 0 over $\mathbb{Q}(\sqrt{2})$. Thus all of its points over
$\mathbb{Q}(\sqrt{2})$ are torsion points, and we can count them by considering the number of
points in various residue fields of $\mathbb{Q}(\sqrt{2})$ (as in~\cite{Si}, VII.3). Note that 7
splits in $\mathbb{Q}(\sqrt{2})$, so $\mathbb{Q}(\sqrt{2})$ has a residue field isomorphic to
$\mathbb{F}_7$. Now $X_0(15)$ has good reduction at the primes above 7, and
$|X_0(15)(\mathbb{F}_7)|=8$. By~\cite{Si}, VII.3.1(b), we see that the size of the torsion group
over $\mathbb{Q}(\sqrt{2})$ divides 8. However, we know that $X_0(15)$ has 8 points over
$\mathbb{Q}$, all of which are torsion, and so these can be the only points on $X_0(15)$ defined
over $\mathbb{Q}(\sqrt{2})$. Of these, 4 are cusps, and the remaining 4 correspond to elliptic
curves over $\mathbb{Q}$ which have additive reduction at 5 (curves 50A1, 50A2, 50A3 and 50A4 in
Cremona's tables). Since 5 is unramified in $\mathbb{Q}(\sqrt{2})/\mathbb{Q}$, these curves
continue to have additive reduction at 5 over $\mathbb{Q}(\sqrt{2})$, and so are also not
semistable. It follows that none of the $\mathbb{Q}(\sqrt{2})$-rational points on $X_0(15)$
correspond to semistable elliptic curves.
\end{Proof}

\begin{theorem}\label{th9.1}
Any semistable elliptic curve over $\mathbb{Q}(\sqrt{2})$ is modular.
\end{theorem}

\begin{Proof}
Let $E$ be a semistable elliptic curve over $\mathbb{Q}(\sqrt{2}).$ By Proposition~\ref{prop9.6},
one of $\overline\rho_{E,3}$ or $\overline\rho_{E,5}$ will be absolutely irreducible. The case
where $\overline\rho_{E,3}$ is absolutely irreducible and $E$ has good ordinary reduction or
multiplicative reduction at $3$ follows from Theorem~5.1 of~\cite{SW2} (using
Proposition~\ref{prop9.5} to check the hypothesis that
$\overline\rho_{E,3}|_{\mathrm{Gal}(\overline{F}/F(\sqrt{-3}))}$ is absolutely irreducible). If
$\overline\rho_{E,3}$ is absolutely irreducible and $E$ has supersingular reduction, then the
modularity of $E$ follows from Proposition~\ref{prop9.3}. Otherwise $\overline\rho_{E,5}$ is
irreducible, and modularity follows by switching to an elliptic curve $E'$ as in the proof of
Proposition~\ref{prop9.4}. By the previous argument, $E'$ is modular, so that
$\overline\rho_{E',5}\cong\overline\rho_{E,5}$ is modular. If $E$ has good ordinary reduction or
multiplicative reduction at $5$, modularity follows from Theorem~5.1 of~\cite{SW2}, again using
Proposition~\ref{prop9.5} to check that the hypotheses of this theorem hold. Otherwise, $E$ has
good supersingular reduction at $5$. As remarked at the end of \S\ref{modularity}, since $5$ is
unramified in $\mathbb{Q}(\sqrt{2})$, the Galois representation $\overline\rho_{E,5}$ has the
form given in Theorem~\ref{th7.1}; this theorem now implies that $E$ is modular, as required.
\end{Proof}

\begin{remark} \rm
In fact, $\mathbb{Q}(\sqrt{2})$ is not the only real quadratic field for which all the numerology
is valid to deduce modularity. Indeed, let $F=\mathbb{Q}(\sqrt{17})$. Note that 3 and 5 are inert
in $F$. Again using {\tt PARI}, one can verify that the class numbers of $F(\zeta_3)$ and
$F(\zeta_5)$ are both 1, so that the analogue of Proposition~\ref{prop9.5} will hold also for
$F$. (We suspect that this might be the only other real quadratic field with this property.)
Next, the quadratic twist of $X_0(15)$ to $F$ is curve 4335D3 in Cremona's tables, which has rank
0 (and 4 points defined over $\mathbb{Q}$), so that $X_0(15)$ has rank 0 over $F$. We can count
the $\mathbb{Q}(\sqrt{17})$-rational points by counting the points in residue fields of $F$ whose
characteristic is a prime of good reduction for $X_0(15)$. Since 13 and 43 both split in $F$, and
$X_0(15)$ has 16 points in $\mathbb{F}_{13}$ and 40 points in $\mathbb{F}_{43}$, we see that the
size of the torsion group of $X_0(15)$ over $F$ divides 8. Now one argues as in the case of
$\mathbb{Q}(\sqrt{2})$ to see that all semistable elliptic curves over $\mathbb{Q}(\sqrt{17})$
are modular. \end{remark}

\appendix
\section{The Honda system associated to a Raynaud scheme}
We describe the Honda system associated to a Raynaud scheme.
We fix:
\begin{eqnarray*}
k& :& \mbox{a perfect field of odd characteristic}\ p, \\
W(k),v& :& \mbox{its Witt ring and normalized valuation}\ v \
(\mathrm{so}\ v(p)=1),\\
\sigma & :& \mbox{the Frobenius automorphism}\
\sigma:W(k)\longrightarrow W(k).
\end{eqnarray*}
We also fix a finite field $\mathbb{F}$ of order $p^r,$ and we
assume that there is an injection $\mathbb{F}\hookrightarrow k$ of
fields. There are $r$ `fundamental characters', indexed by a
principal homogeneous space over $\mathbb{Z}/r\mathbb{Z}.$ We
recall the definition (D\'{e}finition 1.1.1 of \cite{raynaud}):
these are maps $\chi_i:\mathbb{F}\rightarrow W(k)$ such that
$\chi_i|_{\mathbb{F}^\times}$ is a multiplicative character,
$\chi_i(0)=0,$ and the composite
$$\mathbb{F}\stackrel{\chi_i}{\longrightarrow}W(k)\stackrel{\bmod
p}{\longrightarrow}k$$ is a homomorphism of fields.

Given a multiplicative character
$\chi:\mathbb{F}^\times\rightarrow W(k)^\times,$ Raynaud defines a
quantity $\omega_\chi$ in $W(k).$ Raynaud also defines, starting
from a fundamental character, another quantity $\omega\in W(k)$
(which is then shown to be independent of the choice of
fundamental character).  We refer to equations 11 of
\cite{raynaud} for the defining relations, and we recall that
(Proposition 1.3.1 of \cite{raynaud}):
\begin{itemize}
\item Write $\chi$ (uniquely) as a product over fundamental
characters $\prod_{i\in\mathbb{Z}/r\mathbb{Z}}\chi_i^{a_i}$ with
$0\leq a_i\leq p-1.$ Then
$$\omega_\chi\equiv a_1!\ldots a_r!~(\bmod p).$$
\item $\omega\equiv p!~(\bmod p^2).$
\end{itemize}
Given two multiplicative characters
$$\chi'=\prod_{i\in\mathbb{Z}/r\mathbb{Z}}\chi_i^{a_i'}\quad\mathrm{and}
\quad\chi''=\prod_{i\in\mathbb{Z}/r\mathbb{Z}}\chi_i^{a_i''}\quad\mathrm{with}\quad
0\leq a_i',a_i''\leq p-1$$ and such that $\chi'\chi''=\chi_i$ for
some $i\in\mathbb{Z}/r\mathbb{Z},$ we let $h$ be the unique
integer such that $0<h\leq r$ and
\begin{eqnarray*}
a_{i-h}'+a_{i-h}'' &=& p,\\
a_{i-k}'+a_{i-k}'' &=& p-1 \quad \mathrm{for} \quad 0< k<h,\\
a_j'&=&a_j''\quad = 0\quad \mathrm{otherwise}.
\end{eqnarray*}

Throughout we fix a Raynaud scheme $G$ over $\mathrm{Spec}\,W(k),$
and we denote its coordinate ring by $A.$ We also fix the
following presentation of $A:$\ it is generated as a
$W(k)$-algebra by $X_i$ with $i\in\mathbb{Z}/r\mathbb{Z}$ and
relations
$$X_i^p=\delta_iX_{i+1}\ \mathrm{where}\ v(\delta_i)=0\
\mathrm{or}\ 1.$$ We also set $\gamma_j=\omega/\delta_j$ and
denote by $\lambda_j$ the Teichm\"{u}ller lift of $\gamma_j.$ The
comultiplication map $\Delta:A\longrightarrow A\otimes A$ is given
by the formula
$$\Delta(X_i)=X_i\otimes 1+1\otimes
X_i+\sum_{\chi'\chi''=\chi_i}
\frac{\gamma_{i-h}\ldots\gamma_{i-1}}{\omega_{\chi'}\omega_{\chi''}}
\left(\prod_{j}X_j^{a_j'}\right)\otimes\left(\prod_{j}X_j^{a_j''}\right).$$

\begin{theorem}
Let $G_{/\mathrm{Spec}\,W(k)}$ be as above. Let $M$ be the
$r$-dimensional $k$-vector space with basis given by $\mathrm{\bf
e}_i, i\in\mathbb{Z}/r\mathbb{Z}.$ Define Frobenius semi-linear
maps $F, V:M\longrightarrow M$ by setting
\begin{eqnarray*}
F(\mathrm{\bf e}_i)&=&\delta_i\mathrm{\bf e}_{i+1}\quad\mathrm{and}\\
V(\mathrm{\bf e}_i)&=&\lambda_{i-1}^{p^{-1}}\mathrm{\bf e}_{i-1},
\end{eqnarray*}
on the basis elements $\mathrm{\bf e}_i$ and then extending
semi-linearly. (So $F(\alpha v)=\sigma(\alpha)F(v)$ and $V(\alpha
v)=\sigma^{-1}(\alpha)V(v)$ for $\alpha\in k.$) Let $L\subset M$
be the $k$-linear subspace of $M$ spanned by
$\lambda_{i-1}\mathrm{\bf e}_i,i\in\mathbb{Z}/r\mathbb{Z}.$ Then
$(L,M)$ is the Honda system associated to $G.$
\end{theorem}

\subsection*{Witt covectors}

For $n\geq 0,$ set
$$W_n=W_n(X_0,\ldots
,X_n)\stackrel{\mathrm{def}}{=}X_0^{p^n}+pX_1^{p^{n-1}}+\cdots
+p^nX_n.$$ There are polynomials
$$S_n\in \mathbb{Z}[Y_0,Y_1,\ldots ; Z_0,Z_1,\ldots]$$ where
$S_0=Y_0+Z_0$, and $S_n,$ for $n\geq 1,$ satisfies the relation
$$W_n(S_0,S_1,\ldots ,S_n)=W_n(Y_0,Y_1,\ldots ,Y_n)+W_n(Z_0,Z_1,
\ldots, Z_n).$$

We record the following for future use.
\begin{proposition}
 Modulo the ideal $(p,Y_0^p,\ldots , Y_{n-2}^p,Z_0^p,\ldots , Z_{n-2}^p),$ we have
\begin{eqnarray*}
\lefteqn{S_n(Y_0,\ldots ,Y_n;Z_0, \ldots ,Z_n)=
 Y_n+Z_n +\sum_{i=1}^{p-1}
\frac{Y_{n-1}^iZ_{n-1}^{p-i}}{i! (p-i)!}+}\\
 & & \sum_{r=0}^{n-2}(-1)^{n-r}\left((Y_{n-1}+Z_{n-1})\cdots
(Y_{r+1}+Z_{r+1})\right)^{p-1}\sum_{i=1}^{p-1}
\frac{Y_r^iZ_r^{p-i}}{i! (p-i)!}.
\end{eqnarray*}
\end{proposition}

\begin{Proof}
 We have
$$S_n=Y_n+Z_n+\frac{Y_{n-1}^p+Z_{n-1}^p-S_{n-1}^p}{p}
+\frac{Y_{n-2}^{p^2}+Z_{n-2}^{p^2}-S_{n-2}^{p^2}}{p^2}+\ldots $$
Thus modulo $X_i^p,\ i\geq 0,$ we have
$$S_n =Y_n+Z_n+\frac{Y_{n-1}^p+Z_{n-1}^p-S_{n-1}^p}{p}.$$
Assume the proposition for $n-1.$ The right hand side of the above
relation, modulo $(p,X_0^p,X_1^p,\ldots ),$ is equal to
\begin{eqnarray*}
\lefteqn{Y_n+Z_n+\frac{Y_{n-1}^p+Z_{n-1}^p-(Y_{n-1}+Z_{n-1})^p}{p}}\\
& &
-(Y_{n-1}+Z_{n-1})^{p-1}\sum_{r=0}^{n-2}\left\{(-1)^{n-1-r}\prod_{j=r+1}^{n-2}(Y_{j}+Z_{j})^{p-1}
\sum_{i=1}^{p-1} \frac{Y_r^iZ_r^{p-i}}{i!
(p-i)!}\right\}.\end{eqnarray*} This proves the
proposition.\end{Proof}

\begin{definition}\rm
We define $\widetilde{S_{-n}}\in \mathbb{F}_p[Y_{-n},\ldots
,Y_0;Z_{-n},\ldots ,Z_0]$ to be the polynomial
\begin{eqnarray*}
\lefteqn{Y_0+Z_0 +\sum_{i=1}^{p-1} \frac{Y_{-1}^iZ_{-1}^{p-i}}{i!
(p-i)!}+}\\ & &
\sum_{r=2}^{n}(-1)^{r-1}\left((Y_{-1}+Z_{-1})\cdots
(Y_{-r+1}+Z_{-r+1})\right)^{p-1}\sum_{i=1}^{p-1}
\frac{Y_{-r}^iZ_{-r}^{p-i}}{i! (p-i)!}.
\end{eqnarray*}
We have $S_n(Y_{-n},\ldots ,Y_0;Z_{-n},\ldots
,Z_0)\equiv\widetilde{S_{-n}}\bmod (p,Y_{-2}^p,\ldots ,
Y_{-n}^p,Z_{-2}^p,\ldots , Z_{-n}^p).$
\end{definition}

The following
definition is due to Fontaine~(\cite{fontaine}).
\begin{definition}\rm
 For any finite $k$-algebra $R,$ the group of
$R$-valued Witt covectors $CW_k(R)$ is given by:
 \begin{itemize}
\item As a set, the elements of $CW_k(R)$ are sequences
$$\left\{(\ldots , a_{-n},\ldots ,a_{-1},a_0):a_{-i}\in R\ \mbox{is
nilpotent for large}\ i\right\}.$$
\item For $(a_{-i})_{i\geq 0},(b_{-i})_{i\geq 0}\in CW_k(R),$ let
$$c_{-n}=\lim_{m\to\infty}S_m(a_{-n-m},\ldots ,a_{-n};b_{-n-m},\ldots
,b_{-n}).$$ The sequence $(c_{-i})_{i\geq 0}\in CW_k(R),$ and the
group law is
$$(a_{-i})_{i\geq 0}+(b_{-i})_{i\geq 0}\stackrel{\mathrm{def}}{=}(c_{-i})_{i\geq
0}.$$
\item The identity element is $(\ldots ,0,0).$
These give $CW_k(R)$ the structure of a commutative group
(Proposition 1.4, Chapter II of~\cite{fontaine}).
\end{itemize}
\end{definition}

$CW_k(R)$ has a natural structure of a $W(k)$-module which, for
$x\in k,$ is given by
$$[x](\ldots , a_{-n},\ldots ,a_{-1}, a_0)=(\ldots ,
x^{p^{-n}}a_{-n},\ldots ,x^{p^{-1}}a_{-1}, xa_0).$$ Here, $[x]$ is
the Teichm\"{u}ller lift of $x.$ The Frobenius and Verschiebung
operators $F,V:CW_k(R)\longrightarrow CW_k(R)$ are given by
\begin{eqnarray*}
F(\ldots , a_{-n},\ldots ,a_{-1}, a_0)&\stackrel{\mathrm{def}}{=}
&(\ldots , a_{-n}^p,\ldots ,a_{-1}^p, a_0^p), \ \mathrm{and}\\
 V(\ldots ,a_{-n},\ldots ,a_{-1}, a_0)&\stackrel{\mathrm{def}}{=}&
(\ldots ,a_{-n+1},\ldots , a_{-1}).
\end{eqnarray*}
These are additive, and they satisfy the relation $FV=VF=p.$ As
for compatibility with the $W(k)$-module structure, one has
$F\alpha=\sigma (\alpha)F$ and $V\alpha=\sigma^{-1} (\alpha)V$
where $\sigma :W(k)\longrightarrow W(k)$ is the Frobenius. In
other words, the Witt covectors $CW_k(R)$ form a module over the
Dieudonn\'{e} ring $D_k=W(k)[F,V].$

\subsection*{The Dieudonn\'{e} module of the special fibre}

We now calculate the Dieudonn\'{e} module associated to the
special fibre $G_k.$ We know that this is a vector space over $k$ of
dimension $r.$

Following \cite{fontaine}, we need to calculate certain elements of
$A_k$-valued Witt covectors. These elements are formal group
homomorphisms from $G_k$ to $\widehat{CW_k}$.
We will describe the `homomorphism' condition
shortly, but let us start with a candidate covector $\mathrm{\bf
a}= (\ldots ,a_{-n},\ldots ,a_0).$ Since $FV(\mathrm{\bf
a})=p\mathrm{\bf a}=0,$ we must have $a_{-n}^p=0$ for $n\geq 1.$
The comultiplication map
$$\Delta:A_k\longrightarrow A_k\otimes_kA_k$$
gives us an $A_k\otimes_kA_k$-valued Witt covector
$$\Delta(\mathrm{\bf
a})= (\ldots ,\Delta(a_{-n}),\ldots ,\Delta(a_0)).$$ In
$CW_k(A_k\otimes_kA_k),$ we also have the sum
$$\mathrm{\bf
a}\otimes 1+1\otimes\mathrm{\bf a}=(\ldots ,a_{-n}\otimes 1,\ldots
,a_0\otimes 1)+(\ldots ,1\otimes a_{-n},\ldots ,1\otimes a_0).$$
The `homomorphism' condition is then
$$\Delta(\mathrm{\bf
a})=\mathrm{\bf a}\otimes 1+1\otimes\mathrm{\bf a}.$$

We now define covectors
$$\mathrm{\bf
e}_i\stackrel{\mathrm{def}}{=}\left(\ldots,
\underbrace{\gamma_{i-1}^{p^{-n}}\cdots\gamma_{i-n+1}^{p^{-2}}
\gamma_{i-n}^{p^{-1}}}_{n\ \mathrm{factors} }X_{i-n}, \ldots
,\gamma_{i-1}^{p^{-1}}X_{i-1},X_i\right)$$ for
$i\in\mathbb{Z}/r\mathbb{Z}.$ Here we are viewing $X_j$ and
$\gamma_j$ modulo $p$ (and so $\gamma_j^{p^{-n}}$ is the mod $p$
reduction of $\sigma^{-n}(\gamma_j)$). We shall check (by a
tedious, but entirely straightforward, calculation) that
$\mathrm{\bf e}_i$ satisfies the `homomorphism' condition. Note
that $F\mathrm{\bf e}_n=(\ldots, 0,\ldots , 0,X_n^p)$ as
$X_n^p\neq 0\bmod p$ implies $\gamma_n\equiv 0\bmod p$.

We write $Y_i=X_i\otimes 1$ and $Z_i=1\otimes X_i,$ and so
\begin{eqnarray*}
\mathrm{\bf e}_i\otimes 1 &=&
(\ldots,\gamma_{i-1}^{p^{-1}}Y_{i-1},Y_i),\\
1\otimes\mathrm{\bf e}_i
&=&(\ldots,\gamma_{i-1}^{p^{-1}}Z_{i-1},Z_i).\end{eqnarray*}

\begin{lemma}
\begin{itemize}
\item For any $n\geq 0,$ we have
\begin{eqnarray*}
\lefteqn{\gamma_{i-1}\gamma_{i-2}\cdots\gamma_{i-(n+1)}=}\\
& & \left\{\left(\gamma_{i-1}^{p^{-1}}\right)
\left(\gamma_{i-1}^{p^{-2}}\gamma_{i-2}^{p^{-1}}\right)\cdots
\left(\gamma_{i-1}^{p^{-n}}\cdots\gamma_{i-n}^{p^{-1}}\right)\right\}^{p-1}\times
\left(\gamma_{i-1}^{p^{-(n+1)}}\cdots\gamma_{i-(n+1)}^{p^{-1}}\right)^p
\end{eqnarray*}
\item For any $n\geq r-1,$ we have
$$\Delta(X_i)=\widetilde{S_{-n}}(\mathrm{\bf e}_i\otimes 1;1\otimes\mathrm{\bf e}_i).$$
\item $\Delta (\mathrm{\bf e}_i)=\mathrm{\bf e}_i\otimes 1+1\otimes\mathrm{\bf
e}_i.$
\end{itemize}
\end{lemma}

\begin{Proof}
The first part is a simple manipulation of
symbols; the third part follows from the second (for example,
apply the $V$ operator).

We now prove the second part of the lemma. Since $X_k$
(resp.~$Y_k,Z_k,\gamma_k$) is $X_{k+r}$ (resp. $Y_{k+r},
Z_{k+r},\gamma_{k+r}$), it follows that the limit
$$\lim_{m\to\infty}S_m\left(\left(\gamma_{i-1}^{p^{-m}}\cdots
\gamma_{i-m}^{p^{-1}}\right)Y_{i-m}, \ldots ,Y_i;
\left(\gamma_{i-1}^{p^{-m}}\cdots
\gamma_{i-m}^{p^{-1}}\right)Z_{i-m}, \ldots ,Z_i\right)$$ is equal
to
$$\widetilde{S_{-r+1}}(\mathrm{\bf e}_i\otimes 1 ;
1\otimes\mathrm{\bf e}_i)=\widetilde{S_{-n}}(\mathrm{\bf
e}_i\otimes 1 ; 1\otimes\mathrm{\bf e}_i)\ \mathrm{for~any} \
n\geq r-1.$$
Using the first part, we see that $\widetilde{S_{-r+1}}$ is equal to
\begin{eqnarray*}
\lefteqn{Y_i+Z_i+\gamma_{i-1}\sum_{j=1}^{p-1}\frac{Y_{i-1}^jZ_{i-1}^{p-j}}{j!(p-j)!}+}\\
& &
\sum_{k=1}^{r-2}(-1)^{k}\left(\prod_{j=1}^{k+1}\gamma_{i-j}\right)
\left(\prod_{j=1}^{k}(Y_{i-j}+Z_{i-j})^{p-1}\right)
\sum_{j=1}^{p-1}\frac{Y_{i-k-1}^jZ_{i-k-1}^{p-j}}{j!(p-j)!}.
\end{eqnarray*}
Now
\begin{eqnarray*}
\lefteqn{(Y_{i-1}+Z_{i-1})^{p-1}\cdots(Y_{i-k+1}+Z_{i-k+1})^{p-1}
\sum_{j=1}^{p-1}\frac{Y_{i-k}^jZ_{i-k}^{p-j}}{j!(p-j)!}}\\
&=&\sum\frac{(-1)^{k-1}}{a_{i-1}!\cdots a_{i-k}!b_{i-1}!\cdots
b_{i-k}!}Y_{i-1}^{a_{i-1}}\cdots Y_{i-k}^{a_{i-k}}
Z_{i-1}^{b_{i-1}}\cdots Z_{i-k}^{b_{i-k}}
\end{eqnarray*}
where the sum is over
$$\left\{\begin{matrix} 0\leq a_{i-1},\ldots,a_{i-k+1}\leq p-1,\
1\leq a_{i-k}\leq p-1,\quad
\mathrm{and}\\
a_{i-j}+b_{i-j}=p-1,\ j=1,\dots , k-1; \
a_{i-k}+b_{i-k}=p.
\end{matrix}\right.$$
The second part of the lemma easily follows.
\end{Proof}

\begin{lemma} The covectors $\mathrm{\bf e}_1,\ldots, \mathrm{\bf e}_r$
are $k$-linearly independent.
\end{lemma}

\begin{Proof} First suppose that one of the $\gamma_i$ is
divisible by $p.$ Then each $\mathrm{\bf e}_n$ has only finitely
many non-zero terms. Applying the $V$ operator, one reduces the
linear independence of $\mathrm{\bf e}_1,\ldots, \mathrm{\bf e}_1$
to the $k$-linear independence of $X_1,\ldots, X_r$, which is
clear.

In the remaining case, we have
$A_k=k[X_1,\ldots,X_r]/(X_1^p,\ldots,X_r^p).$ The $0$th term of
$\alpha_1\mathrm{\bf e}_1+\cdots +\alpha_r\mathrm{\bf e}_r$ is
$$\alpha_1X_1+\cdots +\alpha_rX_r +\mbox{ an element of }
(X_1,\ldots,X_r)^2,$$ and the lemma follows. \end{Proof}

It now follows that the $k$-linear span of $\mathrm{\bf
e}_1,\ldots, \mathrm{\bf e}_r,$ which is a subspace of the
Dieudonn\'{e} module of $G_k,$ is in fact the whole Dieudonn\'{e}
module (as both are of dimension $r$).This gives the following
proposition:

\begin{proposition} The Dieudonn\'{e} module of $G_k$ is the $r$-dimensional
$k$-vector space
$$k\mathrm{\bf e}_1\oplus\cdots\oplus k\mathrm{\bf e}_r$$
with Frobenius and Verschiebung actions given by
\begin{eqnarray*}
F(\mathrm{\bf e}_i)&=&\delta_i\mathrm{\bf e}_{i+1},\quad\mathrm{and}\\
V(\mathrm{\bf e}_i)&=&\lambda_{i-1}^{p^{-1}}\mathrm{\bf e}_{i-1}
\end{eqnarray*}
on the basis elements $\mathrm{\bf e}_i,\
i\in\mathbb{Z}/r\mathbb{Z},$ which one then extends semi-linearly.
$\hfill\Box$
\end{proposition}

\subsection*{The Honda system associated to $G$}

We need to determine the kernel $L$ of the composite
$$M\hookrightarrow CW_k(A_k)\stackrel{w}{\longrightarrow}A_K/pA$$
where $w$ is defined as follows: for $(a_{-n})\in CW_k(A_k),$
choose for each $a_{-n}$ a lift $\hat{a}_{-n}\in A ,$ and define
$$w((a_{-n}))\stackrel{def}{=}\sum_{n\geq
0}p^{-n}\hat{a}_{-n}^{p^n} \bmod pA.$$ This is well-defined (see
Chapter II, section 5.2 of \cite{fontaine}).

Let $\lambda_j$ be the Teichm\"{u}ller lift of $\gamma_j\bmod p.$
Note that as
$$X_{i-n}^{p^n}=\frac{\omega^{p^{n-1}}}{\gamma_{i-n}^{p^{n-1}}}
\frac{\omega^{p^{n-2}}}{\gamma_{i-n+1}^{p^{n-2}}}\cdots
\frac{\omega}{\gamma_{i-1}}X_i$$
and
$$v\left(\omega^{p^{n-1}+\ldots +1}\right)=(p^n-1)/(p-1),$$
we have
$$p^{-n}\left(\lambda_{i-1}^{p^{-n}}\cdots\lambda_{i-n}^{p^{-1}}\right)^{p^n}X_{i-n}^{p^n}
\in pA$$ for $n\geq 2.$ It follows that
\begin{eqnarray*}
w(\mathrm{\bf e}_i) &\cong& X_i +p^{-1}\lambda_{i-1}X_{i-1}^p\bmod
pA\\
& \cong &\left\{ \begin{matrix} X_i\bmod
pA\quad\mathrm{if}\ \lambda_{i-1}=0\\ 0\bmod pA\quad\mathrm{if}\
\lambda_{i-1}\neq 0\end{matrix}\right. \ .
\end{eqnarray*}

Thus $L$ contains the subspace of $M$ spanned by the covectors
$\lambda_{i-1}\mathrm{\bf e}_i.$ Alternatively, $L$ contains the
subspace generated by $\mathrm{\bf e}_i$ with $p$ dividing
$\delta_{i-1}.$ Since $FM$ is the $k$-span of $\mathrm{\bf e}_i$
with $p$ not dividing $\delta_{i-1},$ a dimension count shows
that $L$ is in fact the $k$-span of $\lambda_{i-1}\mathrm{\bf
e}_i.$

\noindent{\sc Current address:} Department of Pure Mathematics, University of Sheffield,
Sheffield S3 7RH, U.K.

{\tt a.f.jarvis@shef.ac.uk}, {\tt j.manoharmayum@shef.ac.uk}

\end{document}